\documentclass[final,1p,times]{elsarticle}

\journal{submit to xxx }

\usepackage{mathrsfs,amsmath,amssymb}

\usepackage{bm}
\usepackage{hyperref}
\usepackage{graphicx}
\usepackage{subfigure}
\usepackage{caption}
\usepackage{overpic}
\usepackage{epsfig}
\usepackage{xcolor}
\usepackage{indentfirst}
\usepackage{fancyhdr}
\usepackage{comment}
\usepackage{makecell, rotating}
\usepackage{algorithm} 
\usepackage{algorithmic} 
\usepackage{float}
\usepackage{fullpage} 
\usepackage{multirow}
\usepackage{booktabs}
\usepackage{xkeyval}
\usepackage{todonotes}
\usepackage{mathtools}
\usepackage{exscale}
\usepackage{relsize}
\presetkeys{todonotes}{inline}{} 
\mathtoolsset{showonlyrefs}

\usepackage[draft]{changes}
\usepackage{lipsum}
\definechangesauthor[name={Huayi Wei}, color=red]{why}
\definechangesauthor[name={Kai Jiang}, color=blue]{JK}
\definechangesauthor[name={Xin Wang}, color=blue]{wangxin}
\usepackage{exscale}
\usepackage{relsize}
\usepackage{multirow}
\usepackage{diagbox}
\usepackage{tikz,mathpazo}
\usepackage{geometry}
\usetikzlibrary{arrows,shapes,chains}
\newcommand{\bbR}{\mathbb{R}}
\newcommand{\bbZ}{\mathbb{Z}}

\newcommand{\br}{\bm{r}}
\newcommand{\be}{\bm{e}}

\newcommand{\bn}{\bm{n}}
\newcommand{\bal}{\bm{\alpha}}

\newcommand{\bM}{\boldsymbol{M}}
\newcommand{\bA}{\boldsymbol{A}}
\newcommand{\bF}{\boldsymbol{F}}

\newcommand{\calP}{\mathcal{P}}
\newcommand{\calM}{\mathcal{M}}

\newcommand{\bq}{\bm{q}}

\begin{document}

\title{
An adaptive virtual element method for the polymer
self-consistent field theory }

\author{Huayi Wei}
\author{Xin Wang}
\author{Chunyu Chen}
\author{Kai Jiang\corref{cor}}
\address{
 School of Mathematics and Computational Science, 
 Hunan Key Laboratory for Computation and Simulation in Science and Engineering,
 Xiangtan University, Xiangtan, Hunan, P.R. China, 411105.
 }
 \cortext[cor]{Corresponding author. Email: kaijiang@xtu.edu.cn.}

\date{\today}

\begin{abstract}
    In this paper, we develop a high-order adaptive virtual element method (VEM) to
    simulate the self-consistent field theory (SCFT) model in arbitrary domains.
    The VEM is very flexible in handling general polygon elements and can treat
    hanging nodes as polygon vertices without additional processing.  Besides,
    to effectively simulate the phase separation behavior in strong segregation
    systems, an adaptive method on polygonal mesh equipped with a new marking
    strategy is developed. This new marking strategy will indicate the times of
    marked elements to be refined and coarsened, making full use of the
    information contained in the current numerical results. Using the halfedge
    data structure, we can apply the adaptive method to the arbitrary polygonal
    mesh.  Numerical results demonstrate that the developed method is efficient
    in simulating polymers' phase behavior in complex geometric domains. The
    accuracy is consistent with theoretical results.  The adaptive method can
    greatly reduce computational costs to obtain prescribed numerical accuracy
    for strong segregation systems. 
\end{abstract}

\maketitle

\section{Introduction}
\label{sec:intrd}

Block polymers have attracted considerable attention for many years due to their
industrial applications relying on customized microstructures. There are many
industrial applications for the block-copolymer ordered structures at the
nanoscale, such as the construction of high-capacity data storage devices,
waveguides, quantum dot arrays, dielectric mirrors, nanoporous membranes,
nanowires, and interference
lithography\,\cite{charlotte2011interplay,segalman2005patterning}.  In the
practical environment, geometric restriction strongly influences the formation
of microstructures, which also provides a new opportunity to engineer novel
structures. Concretely speaking, the confining geometries and surface
interactions can result in structural frustration, confinement-induced entropy
loss, and lead to novel morphologies that are not obtained in bulk
systems\,\cite{wu2004composite, shi2013self, deng2016chiral}. 

Modeling and numerical simulation provide a practical means to investigate the
phase separation behavior of polymer systems. Fully atomistic and
coarse-graining approaches are both computational intensive methods for
calculating equilibrium microstructures of polymer systems, especially for
larger and more complicated geometries\,\cite{wang2015molecular, sethuraman2014}. A more
and effective continuum approach is the self-consistent field theory (SCFT),
which is one of the most successful modern tools for studying the phase
behaviors of inhomogeneous polymer systems, such as self-assembly and 
thermodynamic stability. SCFT can efficiently describe
polymer architecture, molecular composition, polydispersity, polymer subchain
types, interaction potential, and related information as a series of parameters.
SCFT modeling is started with a coarse-grained chain and microscopic interaction
potentials used in the particle model, then transforms the particle-based model
into a field-theoretic framework, finally derives a mean-field equations system
within saddle-point approximation~\cite{fredrickson2006equilibrium}. 

From the viewpoint of mathematics, the SCFT model is a complicated variational
problem with many challenges, such as saddle-point, nonlinearity,
multi-solutions and multi-parameters. It is not easy to obtain an
analytical solution for this model. A numerical simulation is a feasible tool to
solve SCFT, which usually consists of four parts:
screening initial values\,\cite{xu2013strategy, jiang2010spectral,
	jiang2013discovery}, solving time-dependent partial differential equations
(PDEs)\,\cite{ouaknin2016scft,ceniceros2019efficient}, evaluating (monomer)
density operators\,\cite{ceniceros2019efficient}, and finding
saddle-points\,\cite{ceniceros2004numerical, thompson2004improved,
	jiang2015analytic}. The equilibrium state solution of the SCFT corresponds to an
ordered microstructure. Due to the subtle energy difference among different
ordered patterns in polymer systems, a high order numerical method is strongly
needed.

In the past several decades, spectral methods, especially the Fourier spectral
method, have been the predominant tool for solving the SCFT
model\,\cite{matsen1994stable, rasmussen2002improved, cochran2006stability}.
This approach has high-order precision and is efficient when a spectral
collocation method is found. However, the spectral method uses the global basis
functions to discrete the spatial functions, limiting its applications on the
model defined on complex geometric domains and complex boundary conditions. An
alternative approach uses local basis functions to discretize spatial
functions, such as the finite element method (FEM)\,\cite{ackerman2017finite,
	wei2019finite}. The FEM precision depends on the size and quality of the
mesh and the order of local polynomial basis functions. Combined with the
adaptive method \cite{binev2004}, FEM can obtain a more accurate numerical
solution with less calculation cost. However, there is some inconvenience when
using the FEM with the adaptive method, especially when the adaptive mesh
contains hanging nodes \cite{carstensen2009}, polygonal, or concave mesh
elements.

To address these problems, in this work, we develop an efficient approach to
solve the SCFT model in general domains based on the virtual element method
(VEM). The VEM can be considered as an extension of conforming FEM to polygonal
meshes, which has been developed to solve a variety of PDEs, see \cite{B2013vem,
	L2013vem, lfal2013vem, Antonietti2014vem, Zhao2016vem, Chen2017vem} and
references therein. This paper's contribution contains: (a) formulating the SCFT
problem in real space using a high-order VEM-based variational form, (b)
proposing a new adaptive approach that can make full use of obtained numerical
results, (c) using a halfedge data structure to refine and coarsen general
polygonal grids, (d) the capacity of computing highly segregated systems in
arbitrary areas with the economical computational amount.

The remaining sections are organized as follows. In Sec.\,\ref{sec:scft}, we
give the SCFT model defined in the general domain using the Gaussian diblock
chains as an example. In Sec.\,\ref{sec:method}, we present the high-order
adaptive VEM to solve SCFT in detail.  In Sec.\,\ref{sec:rslt}, we demonstrate
the precision and efficiency of our methods by plenty of numerical
experiments. In Sec.\,\ref{sec:conclusion}, we end with several concluding
remarks and future work.

\section{Self-consistent field theory}
\label{sec:scft}

In this section, we give a brief introduction to the SCFT model for an
incompressible AB diblock copolymer melt on an arbitrary domain $\Omega$. 
We consider a system with $n$ conformationally symmetric diblock copolymers and
each has A and B arms joined together with a covalent bond. The total
degree of polymerization of a diblock copolymer is $N$, the A-monomer
fraction is $f$, and the B-monomer fraction is $1-f$. The field-based
Hamiltonian within mean-field approximation for the incompressible diblock
copolymer melt is\,\cite{fredrickson2006equilibrium, cochran2006stability} 
\begin{align}
H =
\frac{1}{|\Omega|}\int_{\Omega} \left\{-w_+(\br) + \frac{w_-^2(\br)}{\chi
	N}\right\}\,d\br-\log Q[w_+(\br), w_-(\br)],
\label{Eqn:hamilton:2}
\end{align}
where $\chi$ is the Flory-Huggins parameter to describe the interaction between
segments A and B. The terms $w_+(\br)$ and $w_-(\br)$ can be viewed as
fluctuating pressure and exchange chemical potential fields, respectively. The
pressure field enforces the local incompressibility, while the exchange
chemical potential is conjugate to the difference of density operators.  $Q$ is
the single chain partition function, which can be computed according to 
\begin{align}
Q = \frac{1}{|\Omega|}\int_\Omega q(\br,s)q^{\dag}(\br,s)\,d\br,
\quad \forall s\in[0,1].
\end{align}
The forward propagator $q(\br,s)$ represents the probability weight that the
chain of contour length $s$ has its end at position $\br$. The variable $s$ is
used to parameterize each copolymer chain such that $s = 0$ represents the tail
of the A block and $s = f$ is the junction between the A and B blocks. 
According to the flexible Gaussian chain model\,\cite{fredrickson2006equilibrium},
$q(\br, s)$ satisfies the following PDE
\begin{subequations}
	\begin{align}
	\frac{\partial }{\partial s}q(\br,s) &=
	R_g^2\nabla^2_{\br} q(\br,s)-w(\br,s)q(\br,s), \quad \br\in\Omega,
	\label{Eqn:PDE:PDE}
	\\
	w(\br,s) &= \left\{   
	\begin{array}{rl}
	w_A(\br) = w_+(\br)-w_-(\br), & \quad 0\leq s \leq f,
	\\
	w_B(\br) = w_+(\br) +w_-(\br), &  \quad f\leq s \leq 1,
	\end{array}
	\label{Eqn:PDE:w}
	\right.
	\end{align}
	\label{Eqn:PDE}
\end{subequations}
with the initial condition $q(\br, 0)=1$ and $R_g$ being the radius of
gyration. The above PDE is well-defined by possessing an appropriate
boundary condition. In this work, we consider the homogeneous Neumann boundary
condition
\begin{align}
\frac{\partial }{\partial \bm{n}}q(\br,s) = 0, \quad
\br\in\partial\Omega.
\label{Eqn:PDE:bd}
\end{align}
The backward propagator $q^{\dag}(\br,s)$, which represents the probability
weight from $s=1$ to $s=0$, satisfies Eqn.\,(\ref{Eqn:PDE}) only with the
right-hand side of Eqn.\eqref{Eqn:PDE:PDE} multiplied by $-1$. The initial
condition is $q^{\dag}(\br,1)=1$. The normalized segment density operators
$\phi_A(\br)$ and $\phi_B(\br)$ follow from functional derivatives of $Q$ with
respect to $w_A$ and $w_B$ and the familiar factorization property of
propagators
\begin{align}
\label{eq:phiA}
& \phi_A(\br) = -\frac{V}{Q}\frac{\delta Q}{\delta w_A} = 
\frac{1}{Q}\int_0^f q(\br,s)q^{\dag}(\br,s)\,ds,
\\
\label{eq:phiB}
& \phi_B(\br) = -\frac{V}{Q}\frac{\delta Q}{\delta w_B} = 
\frac{1}{Q}\int_f^1 q(\br,s)q^{\dag}(\br,s)\,ds.
\end{align}
The first-order variations of the Hamiltonian with respect to fields $w_+$ and $w_-$
lead to the mean-field equations
\begin{align}
\label{Eqn:modelII:incomp}
&\frac{\delta H}{\delta w_+} = \phi_A(\br) + \phi_B(\br) -1 = 0,
\\
\label{Eqn:modelII:muMinus}
&\frac{\delta H}{\delta w_-} = \frac{2 w_-(\br)}{\chi N} - [\phi_A(\br)-\phi_B(\br)] = 0.
\end{align}

The equilibrium state, i.e., $\delta H/\delta w_\pm = 0$, of the SCFT model
corresponds to the ordered structure. Within the standard framework of SCFT,
finding the stationary states requires the self-consistent iterative procedure, as shown in the following flowchart.
\begin{figure}[H]
	\begin{center}	
		\label{fig:scftiter}
		\tikzstyle{startstop} = [rectangle, rounded corners, minimum width=3cm, minimum height=1cm, text centered, text width=5.5cm, draw=black, fill=yellow!20]
		\tikzstyle{io} = [rectangle, rounded  corners,minimum width=4cm,minimum height=1cm,text centered,draw=black,fill=blue!30]
		\tikzstyle{ioo} = [rectangle, rounded  corners,minimum width=3cm,minimum height=1cm,text centered,draw=black,fill=blue!40]
		\tikzstyle{process1} = [rectangle, rounded  corners, minimum width=3cm,minimum height=1cm,text width=7cm, text centered,draw=black,fill=orange!30]
		\tikzstyle{process2} = [rectangle, rounded  corners, minimum width=3cm,minimum height=1cm,text width=5.5cm, text centered,draw=black,fill=green!25]
		\tikzstyle{decision} = [trapezium, trapezium left angle=80, trapezium right angle=110, rounded corners, minimum width=3cm,minimum height=1cm,text centered,text width=3cm, draw=black,fill=red!30]
		\tikzstyle{arrow} = [thick=50cm,->,>=stealth]
		\begin{tikzpicture}[node distance=1.7cm]
		minimum	\node (start) [startstop] {Given an arbitrary domain $\Omega$ and initial fields $w_+$, $w_{-}$};
		\node (input1) [io,below of=start] {Calculate propagators $q(\br,s)$ and $q^\dag(\br, s)$};
		\node (process1) [process1,below of=input1,yshift=-0.1cm] {Compute $Q$, density operators $\phi_A$ and $\phi_B$, and evaluate the Hamiltonian $H$};
		\node (process2) [process2,below of=process1,yshift=-0.1cm] {Update fields $w_+(\br)$ and $w_-(\br)$};
		\node (decision) [decision,right of =process2,xshift=5cm] {Is Hamilton difference less than a prescribed tolerance ?};
		\node (out1) [ioo,below of=process2] {Converged result};
		\draw [arrow] (start) -- (input1);
		\draw [arrow] (input1) -- (process1);
		\draw [arrow] (process1) -- (process2);
		\draw [arrow,draw=red,thick=3cm] (decision.east) to [out=-400,in=-50] node [right=1cm]{yes} (out1.east);
		\draw [arrow,draw=red] (decision.east) to [ out=400,in=50]node[right=0.1cm]{no} (input1.east);
		\draw [arrow] (process2) -- (decision);
		\end{tikzpicture}
	\end{center}
	\vspace{-1.0cm}
	\caption{Flowchart of SCFT iteration.}
\end{figure}
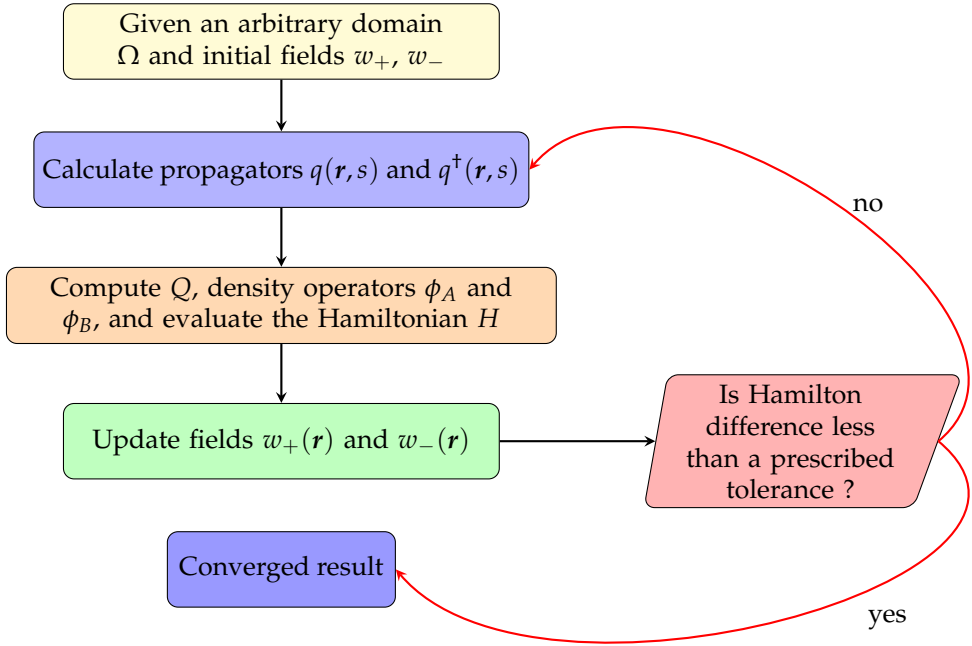
The propagator equation is dependent on the potential fields $w_+(\br)$ and
$w_-(\br)$. In order to start the process, the values of $w_+(\br)$ and
$w_-(\br)$ must be initialized. If the initial values are homogeneous, the gradient
term in the modified diffusion equation goes to zero, leaving no driving force
for forming a microstructure. To prevent this, there must be some
spatial inhomogeneity in the initial values. For a targeted periodic
structure, using the space group symmetry is a useful strategy to screen the
initial configuration\,\cite{jiang2010spectral, jiang2013discovery}. Once
initial values are ready, high-accuracy numerical methods to solve the propagator
equation, and evaluate the density functions, are required to solve the SCFT
model, which is also the main work in this paper. We will detail our approach
in Sec.\,\ref{sec:method}.

The iteration method to update fields is dependent on the mathematical
structure of SCFT. An important fact is that the effective Hamiltonian
\eqref{Eqn:hamilton:2} of diblock copolymers can reach its local minima along
the exchange chemical field $w_-(\br)$, and achieve the maxima along the
pressure field $w_+(\br)$\,\cite{fredrickson2006equilibrium}. Thus alternative
direction gradient approaches, such as the explicit Euler method, can be used
to find the saddle point. In particular, the explicit Euler approach is
expressed as 
\begin{equation}
\begin{aligned}
w_+^{k+1}(\br) &= w_+^k(\br) + \lambda_+
\Big(\phi_A^k(\br)+\phi_B^k(\br)-1 \Big),
\\
w_-^{k+1}(\br) &= w_-^k(\br) - \lambda_-
\left(\frac{2 w_-^k(\br)}{\chi N} - [\phi_A(\br) - \phi_B(\br)]
\right).
\end{aligned}
\end{equation}
An accelerated semi-implicit scheme has been developed to find the equilibrium
states\,\cite{ceniceros2004numerical,jiang2015analytic}. However, the existing
semi-implicit method is based on the asymptotic expansion and global Fourier
transformation and can not be straightforwardly applied to the local basis
discretization schemes.

\section{Numerical methods}
\label{sec:method}
Solving the propagator equations is the most time-consuming part of the entire
numerical simulation, and we will discuss its spatial variables discretization
with the (adaptive) VEM in detail in this section. 
In the following, we use the $\|\cdot\|_{B}$ to denote the common $L^2$ norm over a
finite domain $B$.

\subsection{VEM discretization for the spatial variable}
\label{subsec:VEM}

\begin{figure}[!hbpt]
	\centering
	\includegraphics[width=0.35\linewidth]{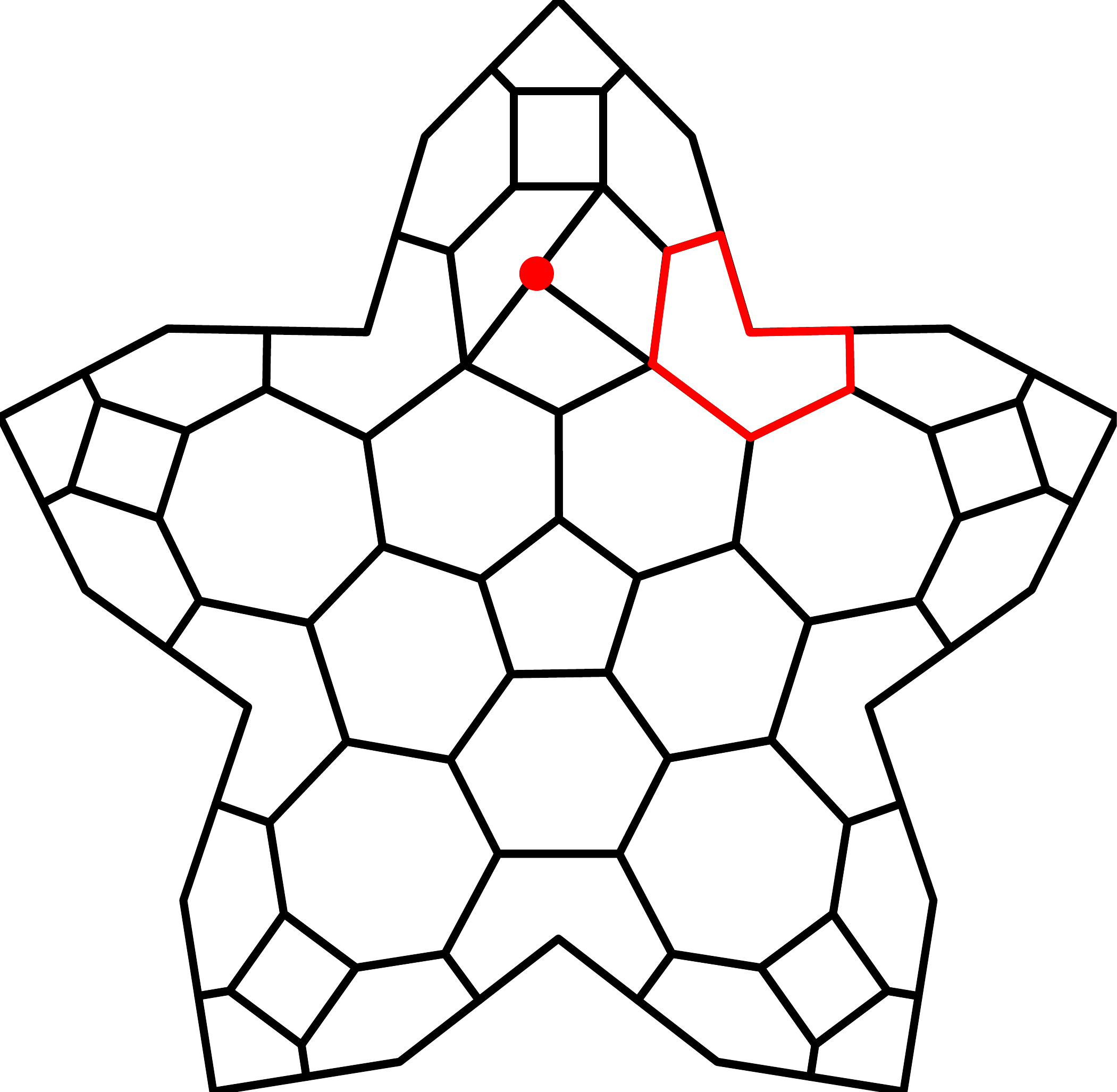}
	\caption{\label{fig:example} A schematic mesh of the VEM including hanging
		node and concave polygons.}
\end{figure}

VEM is a generalization of the finite element method inspired by
the modern mimetic finite difference scheme\,\cite{L2013vem}.
Compared with FEM, VEM can handle general (even non-convex)
polygonal elements. Furthermore, VEM can
naturally treat the handing nodes appearing in the mesh adaptive process as the
vertices of the polygonal elements, which greatly simplifies the design and
implementation of mesh adaptive algorithms. Fig.\,\ref{fig:example} gives a schematic mesh which the VEM can deal with.
Subsequently, we will introduce the virtual element space and discretize 
propagator equations \eqref{Eqn:PDE} based on the variational formulation.

\subsubsection{Virtual element space}
\label{subsubsec:vem}

Let $\Omega_h$ be the polygonal decomposition of a given domain $\Omega \subset
\bbR^2$ including a finite number of non-overlapping polygons. For any polygonal
element $E \in \Omega_h$, let $\partial E=\{e\}$ be the set of all boundary
edges of $E$, $\br_E$ the barycenter, $h_E$ 
the diameter, and $|E|$ the area of the element $E$. Let $\calP_k(E)$ be the 
polynomials space of degree up to $k$ on $E$, $n_k=\dim\calP_k(E) = (n_k+1)(n_k+2)/2$, and
$\calM_k(E) : = \{m_{\bal} : 0 \leq |\bal| \leq k\}$ be the scaled monomial basis set of
$\calP_k(E)$ with form \cite{lfal2013vem}
\begin{equation}
m_{\bal} : = \left(\dfrac{\br - \br_E}{h_E}\right)^{\bal} 
=
\dfrac{(r_1-{r_1}_{E})^{\alpha_1}(r_2-{r_2}_{E})^{\alpha_2}}{{h_E}^{\alpha_1+\alpha_2}},
\quad \alpha_1,~ \alpha_2\in\bbZ_0^+,
\end{equation}
and $|\bal|=\alpha_1+\alpha_2$. We will use $m_\alpha$ instead of
$m_{\bal}$, where $\alpha$ is a one-dimensional index of the natural correspondence of
$\bal$, for example,
\begin{equation}
(0,0) \leftrightarrow 1,\ (1,0) \leftrightarrow 2,\ (0,1) \leftrightarrow 3,\ (2,0)
\leftrightarrow 4,\ \ldots
\end{equation}
The local virtual element space can be defined as~\cite{B2013vem, lfal2013vem}
\begin{equation}
V_{h,E} :=\{v \in H^1(E):~\Delta v \in \calP_{k-2}(E) \
\text{in} \  E; ~ v|_{e} \in \calP_k(e), \forall e \in  \partial E\},
\end{equation}
where $\Delta$ denotes the common Laplace operator. 
$\calP_k(e)$ is a set of polynomials of degree up to $k$ on $e$. The dimension
of $V_{h,E}$ is 
\begin{equation}
N_{dof} = \dim V_{h,E} = n_V +n_V(k-1)+n_{k-2}
\end{equation}
where $n_V$ is the number of vertices of $E$.
The function $v_h \in V_{h,E}$ can be defined by satisfying the following three conditions:
\begin{itemize}
	\item $v_h|_e \in \calP_k(e)$ is a polynomial of degree $ k$ on each edge $e$;
	\item $v_h|_{\partial E} \in C(\partial E)$ is globally continuous
	on $\partial E$;
	\item $\Delta v_h\in \calP_{k-2}(E)$ is a polynomial of degree $ k-2$ in $E$.
\end{itemize}
Correspondingly, the degree of freedom of the $V_{h,E}$ contains:
\begin{itemize}
	\item the value of $v_h$ at the vertices of $E$;
	\item the value of $v_h$ at the $k-1$ internal Gauss-Lobatto
	quadrature points on e;
	\item the moments up to order $k-2$ of $v_h$ in $E$:
	$\frac{1}{E}\int_E v_h m_{\alpha}\,d\br, \alpha = 1,\cdots,n_k$.
\end{itemize} 
Then the global virtual element space can be defined based on the
local space $V_{h,E}$,
\begin{equation}
V_{h} =\{v \in H^1(E):~v|_E \in V_{h,E}, ~\mbox{for all}~ E \in \Omega_h\}.
\end{equation}
The dimension of $V_h$ is 
\begin{equation}
N = \dim{V_h} = N_V + N_E(k-1)+N_P n_{k-2},
\end{equation}
where $N_V$, $N_E $ and $N_P$ are the total number of 
vertices, edges, and elements of $\Omega_h$, respectively.
Since $H^1(V_h)$ is a separable Hilbert space, it can give a set of basis functions
$\{\varphi_i(\br)\}_{i=1}^N$ for $V_h$ such that, for each $u_h(\br)\in V_h$
\begin{align}
u_h(\br) = \sum_{i=1}^N u_i \varphi_i(\br),
\label{eq:VhBasisExp}
\end{align}
where $u_i$ is the coefficient of the degree of freedom corresponding to
$\varphi_i(\br)$. It should be emphasized that the basis functions
$\varphi_i(\br)$ in the VEM do not have explicit expression as the FEM has. In
practical implementation, the quantities related to the basis functions can be
obtained through the degree of freedom.

\subsubsection{Variational formulation}
\label{subsubsec:weak}

Using VEM to solve PDEs \eqref{Eqn:PDE} is based on the variational
formulation whose continuous version is: find $q(\br,s)\in H^1(\Omega)$
such that, for all $v(\br) \in H^1(\Omega)$,
\begin{equation}
\left(\frac{\partial}{\partial s}q(\br,s), v(\br)\right) =
-(\nabla q(\br,s) ,\nabla v(\br)) - (w(\br,s) q(\br,s),
v(\br)), 
\label{eq:continWF}
\end{equation}
where $(\cdot , \cdot)$ represents the $L^2(\Omega)$ inner product.
In numerical computation, the spatial function must be
discretized in the finite-dimensional virtual element space $V_h$. 
Then the continuous variational formulation \eqref{eq:continWF} is discretized
as: find $q_h(\br,s) \in V_h$ such that
\begin{equation}
\label{eqn:discretization}
\left(\frac{\partial}{\partial s}q_h(\br,s), v_h(\br)\right) =
-(\nabla  q_h(\br,s) ,\nabla v_h(\br)) - (w(\br,s) q_h(\br,s) ,
v_h(\br)) , ~  \text{for all} \  v_h(\br) \in V_h.
\end{equation}
Let $v_h(\br) = \varphi_j(\br)$, using the expression \eqref{eq:VhBasisExp}, 
$q_h(\br,s) = \sum_{i=1}^{N} q_i(s) \varphi_i(\br)$.
The discretized variational formulation \eqref{eqn:discretization}
has the matrix form 
\begin{equation}
\label{eq:matrix:form}
\bM \frac{\partial}{\partial s}\bq(s) = -(\bA+\bF)\bq(s),
\end{equation}
where
\begin{equation*}
\bq(s) = (q_1(s), q_2(s), \cdots, q_{N}(s))^T,
\end{equation*}
and
\begin{equation}
\bM_{ij} = (\varphi_i,\varphi_j),\ \bA_{ij} = (\nabla
\varphi_i, \nabla \varphi_j), \ \bF_{ij} = (w(\br,s) \varphi_i,\varphi_j),
\ i,j = 1,\cdots, N.
\end{equation}
The stiffness matrix $\bA$, the mass matrix $\bM$, and the cross mass matrix
$\bF$ can be obtained through projecting local virtual element space $V_{h,E}$
onto polynomial space. In the sequential subsections, we will present the
construction methods for local stiffness, mass, and cross mass matrices. The
corresponding global matrices $\bA$, $\bM$ and $\bF$ can be obtained as the
standard assembly process of FEM once we have the local ones.

\subsubsection{Stiffness matrix}
\label{subsubsec:stiff}

The stiffness matrix in the VEM can be computed by the local $H^1$ projection operator $\Pi^\nabla$, 
\begin{equation}\label{H1proj}
\Pi^\nabla :V_{h,E} \rightarrow \calP_k(E),
\end{equation}
which projects the local virtual element space $V_{h,E}$ onto
the polynomial space with degree up to $k$.  For each $v_h \in V_{h,E}$, we
have the orthogonality condition
\begin{equation}
\label{Eqn:projection}
(\nabla p , \nabla (\Pi^\nabla v_h -v_h))=0, ~\text{for all}~ p \in \calP_k(E).
\end{equation}  
The above condition defines $\Pi^\nabla v_h$ only up to a
constant. It can be fixed by prescribing a projection operator
onto constants $P_0$ requiring 
\begin{equation}
\label{Eqn:P0}
P_0 (\Pi^\nabla v_h -v_h) =0.
\end{equation}
$P_0$ can be chosen as
\begin{equation}
\begin{aligned}
P_0 v_h : &= \frac{1}{n_V}\sum_{i=1}^{n_V} v_h(\br_i), ~\text{when} \ k=1,\\
P_0 v_h : &= \frac{1}{|E|}\int_E v_h
\,d\br=\frac{1}{|E|}(1,v_h)_E, ~\text{when} \quad  k \geq
2,\\
\end{aligned}
\end{equation}
where $n_V$ is the number of vertices of $E$. 

Next we compute the local stiffness matrix $(\bA^{E})_{ij}$ on the polygon $E$, 
\begin{equation}
\label{Eqn:stiffnessmatrix}
(\bA^E)_{ij} = (\nabla \varphi_i, \nabla \varphi_j), ~i,j = 1,\cdots,
N_{dof}.
\end{equation}
With the operator $\Pi^{\nabla}$, $\varphi_i$ can be split into
\begin{equation}
\varphi_i = \Pi ^{\nabla} \varphi_i + (\bm{I} - \Pi ^{\nabla}) \varphi_i,
\end{equation}
Eqn.\,\eqref{Eqn:stiffnessmatrix} becomes
\begin{equation}
(\bA^E)_{ij} = (\nabla \Pi ^{\nabla}\varphi_i, \nabla \Pi^{\nabla}\varphi_j)
+(\nabla ({\bm{I}}-\Pi ^{\nabla})\varphi_i, \nabla ({\bm{I}}-\Pi^{\nabla})\varphi_j).
\end{equation}
Replacing the second term as
$$
S^{E}_0\left(\big({\bm{I}}-\Pi^{\nabla}\big)
\varphi_i,\big({\bm{I}}-\Pi^{\nabla}\big) \varphi_j\right):=\sum_{r=1}^{N_{dof}}
\operatorname{dof}_{r}\left(\big({\bm{I}}-\Pi^{\nabla}\big)\varphi_i\right)
\operatorname{dof}_{r}\left(\big({\bm{I}}-\Pi^{\nabla}\big) \varphi_j\right),
$$
where $\operatorname{dof}_{r}(\varphi_i)=\delta_{ri}$, we obtain the approximate local stiffness matrix
\begin{equation}
(\bA_h^E)_{ij} := (\nabla \Pi ^{\nabla}\varphi_i, \nabla \Pi^{\nabla}\varphi_j)
+S^{E}_0\left(\big({\bm{I}}-\Pi^{\nabla}\big)
\varphi_i,\big({\bm{I}}-\Pi^{\nabla}\big) \varphi_j\right).
\end{equation}
\subsubsection{Mass matrix}
\label{subsubsec:mass}

The mass matrix in the VEM can be obtained from the local $L^2$ projection $\Pi:
V_{h,E}\rightarrow \calP_k(E)$. For each $v_h\in V_{h,E}$,
\begin{equation}
(\Pi v_h,p_k) = (v_h,p_k), ~\forall p_k\in\calP_k(E).
\end{equation}
where $(v_h,p_k)$ can not be calculated directly.
Next, we show how to compute the local mass matrix $\bM^E$ \cite{L2013vem}
\begin{equation}
(\bM^E)_{ij} = (\varphi_i, \varphi_j), \qquad i,j = 1,\cdots,N_{dof}.
\end{equation} 
Similar to the construction method of the stiffness matrix, we define the basis
function $\varphi_i$ through $L^2$ projection operator $\Pi$ 
\begin{equation}
\varphi_i = \Pi \varphi_i + (\bm{I}-\Pi)\varphi_i.
\end{equation} 
Then
\begin{equation}
(\bM^E)_{ij} = (\Pi \varphi_i, \Pi \varphi_j)+((\mathrm{I}-\Pi)\varphi_i, (\mathrm{I} - \Pi)\varphi_j).
\end{equation}
Replacing the second term in the above equation as
$$
S^{E}_1\left(\left({\bm{I}}-\Pi\right)
\varphi_i,\left({\bm{I}}-\Pi\right) \varphi_j\right):=
|E|\sum_{r=1}^{N_{dof}}
\operatorname{dof}_{r}\left(\left({\bm{I}}-\Pi\right)\varphi_i\right)
\operatorname{dof}_{r}\left(\left({\bm{I}}-\Pi\right) \varphi_j\right)
$$
the local mass matrix is approximated as
\begin{equation}
(\bM_h^E)_{ij} := (\Pi\varphi_i,\Pi\varphi_j)
+S^{E}_1\left(\left({\bm{I}}-\Pi\right)
\varphi_i,\left({\bm{I}}-\Pi\right) \varphi_j\right)
\end{equation}

\subsubsection{Cross mass matrix }
\label{subsubsec:cross}

The local cross mass matrix $\bF_{ij}$ on $E$ can be defined as 
\begin{equation}
(\bF_E)_{ij} = (w\varphi_i, \varphi_j), \qquad i,j = 1,\cdots,N_{dof}.
\end{equation} 
Applying the $L^2$ projection $\Pi: V_{h,E}\rightarrow \calP_k(E)$, as defined in the
above Sec.\,\ref{subsubsec:mass}, into the cross term, the local mass matrix can be
calculated as
\begin{equation}
(\bF_E)_{ij} := (\Pi w\Pi \varphi_i, \Pi \varphi_j).
\end{equation}

\subsubsection{Spatial integral }
\label{subsubsec:spaceint}

Here we present the integration approach over
an arbitrary polygon $E$. We divide the polygon $E$ into
triangles $\tau$ by linking two endpoints of each edge and the
barycenter. Then we apply the common Gaussian quadrature in each
triangle, and summarize these integration values. 
\begin{equation}
\begin{aligned}
\int _E f(\br)\,d\br &=\sum_{\tau} \int _{\tau} f(\br) \,d\br
\approx |E|\sum_{\tau}\sum_{j} w_{\tau, j} f(\br_{\tau,j}),
\end{aligned} 
\end{equation}
where $\{\br_{\tau,j}\}$ is the set of quadrature points of
$\tau$, and $\{w_{\tau,j}\}$ the corresponding quadrature weights.

\subsection{Adaptive technique}
\label{subsec:AVEM}

The adaptive method is an important technique to improve the solution's accuracy and reduce computational complexity. The following is the adaptive
process used in SCFT calculation:
\begin{description}
	\item[Step 1] Solve the SCFT model and obtain the numerical solution on the current mesh.
	\item[Step 2] Estimate error on each element from current numerical results.
	\item[Step 3] Mark mesh elements according to the error estimate.
	\item[Step 4] Refine or coarsen the marked elements.
\end{description}
Next, we present some implementation details of the above adaptive process. 

\begin{figure}[H]
	\centering
	\includegraphics[width=1.0\linewidth]{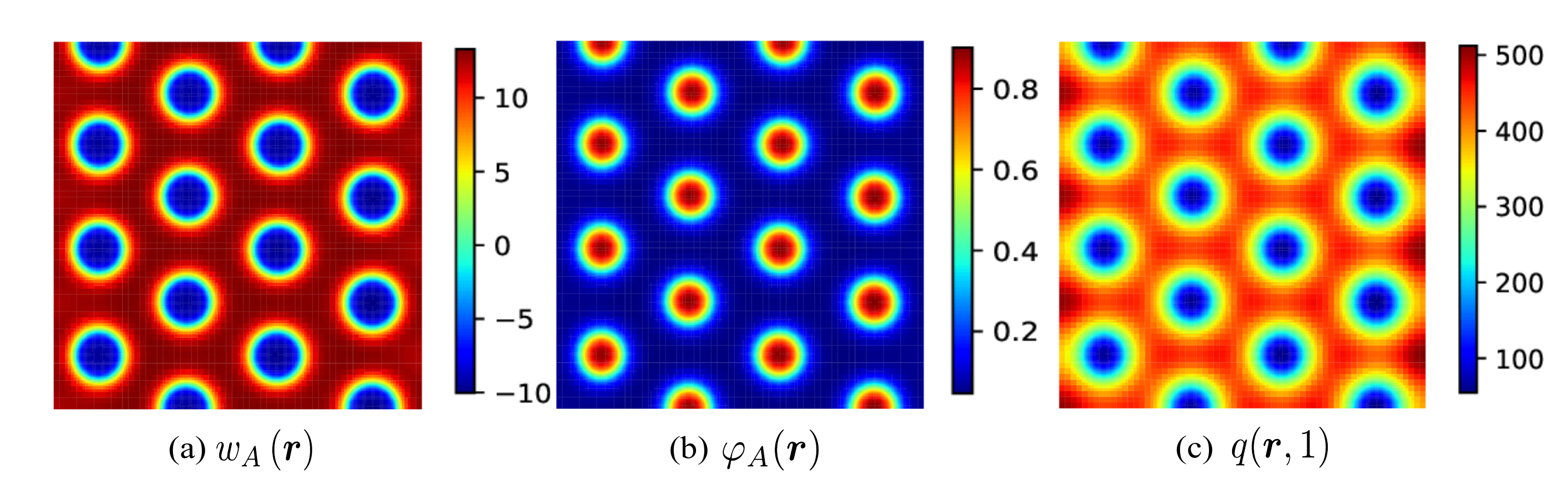}
	\caption{\label{fig:estimator} The equilibrium distributions of
		$w_A(\br)$, $\varphi_A(\br)$, and $q(\br,1)$ when $\chi N =25, f=0.2$.}
\end{figure}

The estimator is an important part of the adaptive method. 
Let $\eta_E$ be the error of the indicator function $u_h$ over each element $E$,
\begin{equation}
\eta_{E} = \|R_h \Pi^\nabla u_h\|_{E},
\label{eq:avem:error}
\end{equation}
$R_h u_h$ is the harmonic average operator \cite{huang2012}
\begin{equation}
R_h u_h := \frac{1}{\sum_{j=1}^{m_z}1/|
	\tau_j|}\sum_{j=1}^{m_z} \frac{1}{|\tau_j|}\nabla \Pi^\nabla u_h \Big|_{\tau_j}.
\end{equation}
$m_z$ is the number of elements $\tau_j$ with $z$ as a vertex.
The indicator function is an essential part of adaptive methods. 
In the SCFT model, several spatial functions can be used as indicator functions,
such as field functions, density functions, and propagators. 
To choose an efficient indicator function, we observe the distribution of these
spatial functions when the SCFT calculation converges.
As an example, Fig.\,\ref{fig:estimator} presents the equilibrium states of $w(\br)$,
$\varphi_A(\br)$ and propagator function of the last contour point $q(\br, 1)$,
respectively, with $\chi N =25, f=0.2$. As one can see, the distributions of
three spatial functions are similar, however, $q(\br,1)$ has the sharpest
interface. If the numerical error of $q(\br,1)$ can be reduced through the
adaptive method, the error of other spatial functions obviously reduces with it.  
Therefore, in the current adaptive method, we choose $q(\br, 1)$ as the
indicator function in the posterior error estimator.

Given an effective and reliable posterior error estimator $\eta_E$, a marking strategy is required to mark mesh elements. Classical marking strategies such as the maximum \cite{jarausch1986} and the $L^2$ criterion
\cite{dorfler1996}, usually refine or coarsen marked mesh elements one time in one adaptive process.
It may make less use of the information of the posterior error estimator.
To improve it, we propose a new marking strategy, named $Log$ criterion, 
as following 
\begin{equation}
n_E = \left[\log_2\frac{\eta_E}{\theta \bar{\eta}}\right],
\end{equation} 
where $\theta$ is a positive constant, $\bar{\eta}$ is the mean value of all
element estimator $\eta_E$, and $[\cdot]$ is the nearest integer function.
$n_E = 0$, $n_E > 0$ and $n_E < 0$ represent that cell $E$ is unchanged, refined
$n_E$ times, and coarsened $|n_E|$ times, respectively.
Obviously, this new $Log$ marking criterion not only denotes which mesh element 
$E$ needs to be improved but also provides the times of refinement or coarseness. 

\begin{figure}[H]
	\centering
	\includegraphics[width=0.5\linewidth]{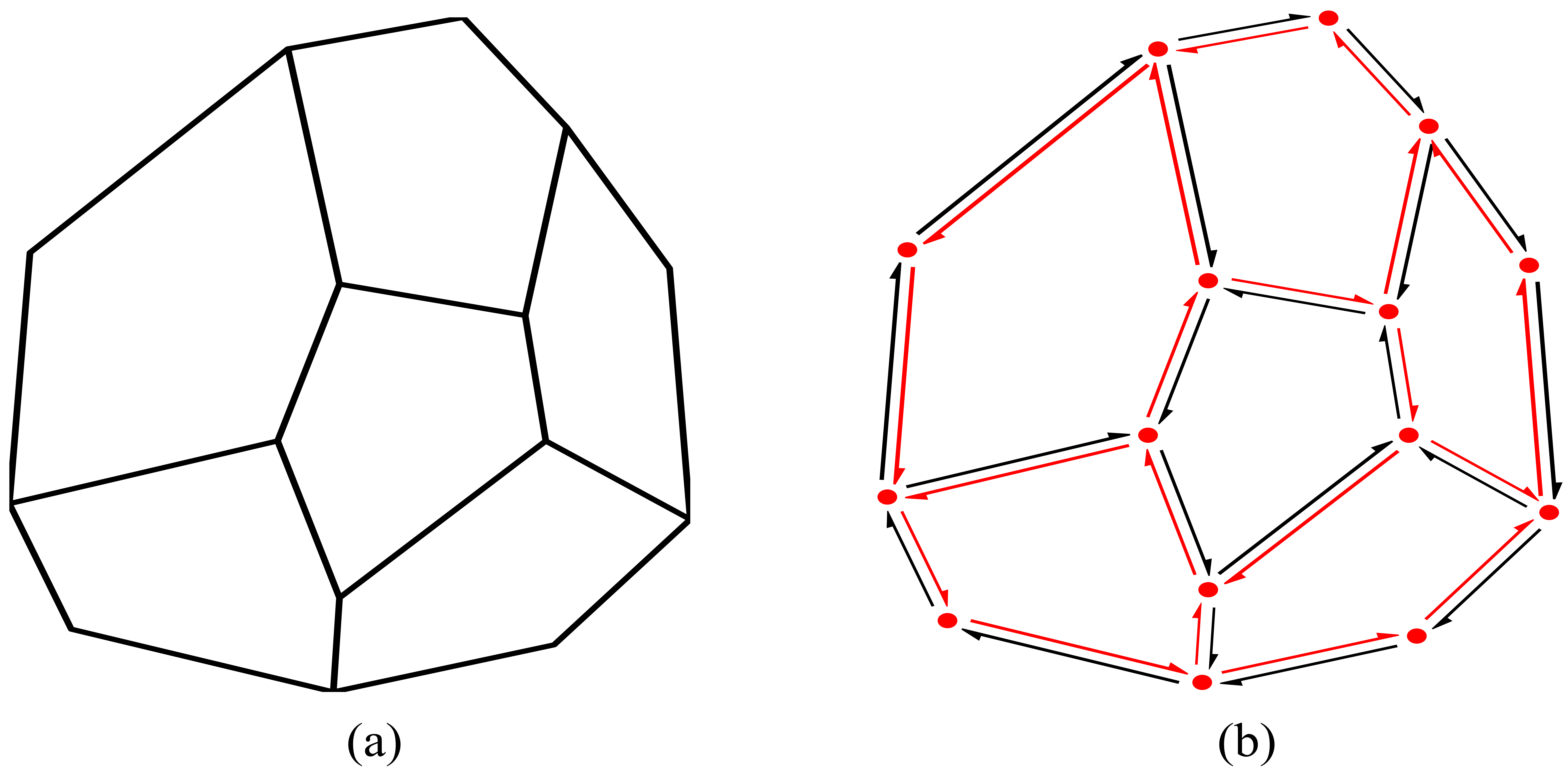}
	\caption{(a) A polygonal mesh. (b) The corresponding
		halfedge data structure of (a).} \label{fig:he}
\end{figure}
\begin{figure}[H]
	\centering
	\includegraphics[width=0.8\linewidth]{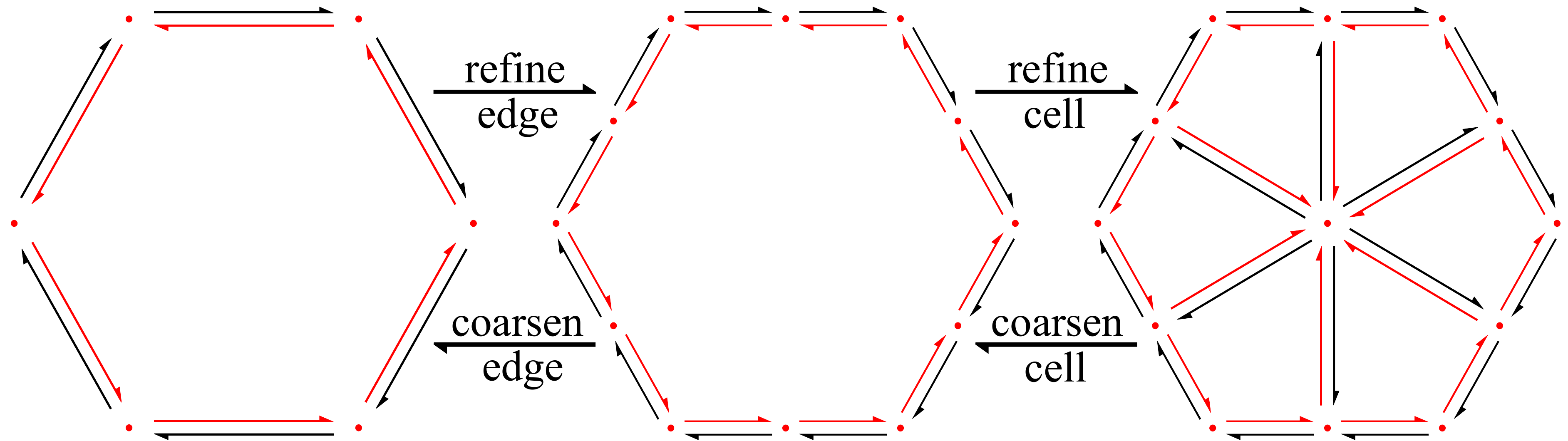}
	\caption{The refinement and coarsening of halfedge mesh.}
	\label{fig:refine}	
\end{figure}

We use the halfedge data structure to implement our adaptive technique which allows us to refine and coarsen arbitrary polygonal mesh.
Halfedge data structure is an edge-centered data structure capable
of maintaining incidence information of nodes, edges, and cells\,\cite{dan2011halfedge}. Each edge is decomposed into two halfedges with opposite orientations. One incident cell and
one incident node are stored in each halfedge. For each cell and each node, one
incident halfedge is stored, see Fig.\,\ref{fig:he}. Halfedge data structure is more
flexible and powerful than the cell-centered data structure. The cell-centered data structure, as a classical data structure, stores the coordinates and indexes of each cell node and requires additional work to reconstruct the relationships between nodes, edges, and cells.
Based on the halfedge data structure, the mesh adaptation is a process of
increasing or decreasing the halfedges, as shown in Fig.\,\ref{fig:refine}.
Notice that, based on the halfedge data structure, the current
mesh refinement and coarsen algorithm, including the red-green approach
\cite{bank1983},  newest vertex bisection \cite{rivara1984} and coarsening
\cite{Chen2010}, can be implemented in a unified way. One can find the
implementation in package FEALPy \cite{fealpy2021}


Once one has the adaptive mesh, the construction of stiffness, mass, cross mass matrices, and the spatial integral formula are the same as Secs.\,\ref{subsubsec:stiff}-\ref{subsubsec:spaceint} present.

\subsection{Full-discrete form}
\label{sec:sdc}

The above matrix form \eqref{eq:matrix:form} is still continuous in the contour variable $s$. There are numerous ways to discretize the contour, such as the second-order operator splitting method, the backward differentiation formulas, Crank-Nicolson (CN) scheme\,\cite{rasmussen2002improved, cochran2006stability, wei2019finite}. Recently, Ceniceros introduced the spectral deferred correction (SDC) approach to improve the accuracy and efficiency of solving polymer SCFT\,\cite{ceniceros2019efficient}. 
In the VEM framework, we choose the SDC scheme to discretize the contour variable. The SDC scheme stemmed from Dutt et al.'s work in 2000\,\cite{Dutt2000sdc},
first solves the PDE with an appropriate method, then uses the residual equation to improve the approximation order of numerical solution.  The key idea of SDC is to use spectral quadratures, such as a Gaussian or a Chebyshev-node interpolatory quadrature, to integrate the contour derivative, which can achieve
a high-accuracy numerical solution with a vastly reduced number of quadrature
points. The detail will be presented in sequential content.

We use the variable step Crank-Nicholson (CN) scheme to solve the
semi-discrete propagator equation \eqref{eq:matrix:form} and obtain the
initial numerical solution $\bq^{[0]}(s)$.  
\begin{equation}
\label{Eqn:CN}
\bM\frac{\bq^{n+1}-\bq^n}{\delta s_n} =
-\frac{1}{2}(\bA+\bF)(\bq^{n+1}+\bq^n),
\end{equation}
where $\delta s_n=s_{n+1}-s_{n}$ is the time step size, $s_n$
($n=0,1,\dots,N_s-1$) is the Chebyshev node\,\cite{clen1960integral}.  It should
be pointed out that other stable time schemes can be employed to solve
semi-discrete propagator equation \eqref{eq:matrix:form}, such as the second-order
operator-splitting method\,\cite{rasmussen2002improved}, implicit-explicit
Runge-Kutta scheme\,\cite{ascher1997rungekutta, ceniceros2019efficient}. 

Then we use the deferred correction scheme to achieve a high-accuracy numerical
solution. We can give the exact semi-discrete solution of propagator by
integrating \eqref{eq:matrix:form} along the contour variable $s$
\begin{equation}\label{Eqn:integral}
\bM\bq(s)=\bM \bq(0)+\int_{0}^{s}\left[(-\bA- \bF) \bq(\tau)\right] d \tau.
\end{equation}
The error between the numerical solution $\bq^{[0]}(s)$ and the exact
semi-discrete solution $\bq(s)$ is defined as
\begin{align}
\be^{[0]}(s) = \bq(s) - \bq^{[0]}(s).
\end{align}
Multiplying both sides by $\bM$, we have
\begin{align}
\bM\be^{[0]}(s) = & \bM\bq(s) - \bM\bq^{[0]}(s) \\
= & \bM \bq(0)+\int_{0}^{s}\left[(-\bA- \bF) \bq(\tau)\right] d \tau -
\bM\bq^{[0]}(s) \\
= & \bM \bq(0)+\int_{0}^{s}\left[(-\bA- \bF) \be^{[0]}(\tau)\right] d \tau + 
\int_{0}^{s}\left[(-\bA- \bF) \bq^{[0]}\right] d \tau- \bM\bq^{[0]}(\tau)(s) \\
= & \int_{0}^{s}\left[(-\bA- \bF) \be^{[0]}\right] d \tau +
\boldsymbol{\gamma}^{[0]}(s),
\end{align}
the residual 
\begin{equation}
\boldsymbol{\gamma}^{[0]}(s) =
\boldsymbol{M}\boldsymbol{q}(0)+\int_{0}^{s}
\left[\left(-\boldsymbol{A}
-\boldsymbol{F}\right)\boldsymbol{q}^{[0]}(\tau)\right]
d\tau-\boldsymbol{M}\boldsymbol{q}^{[0]}(s),
\label{eq:residual}
\end{equation}
can be computed by 
the spectral integral method with Chebyshev-nodes as presented in the Appendix.
By the definition of residual $\bm{\gamma}^{[0]}$, we have the error integration
equation
\begin{align}
\bM\bm{\be}^{[0]}(s) =  \int_0^s
(-\bA-\bF)\bm{\be}^{[0]}(\tau) \,d\tau
+ \bm{\gamma}^{[0]}(s) .
\label{}
\end{align}
Taking the first derivative of the above equation with respect to
$s$ leads to 
\begin{align}
\bM\frac{d{\bm\be}^{[0]}}{ds} = (-\bA-\bF){\bm\be}^{[0]}(s)
+\frac{d {\bm\gamma}^{[0]}}{ds},
\label{eq:erreq}
\end{align}
which can also be solved by the variable step CN scheme\,\eqref{Eqn:CN}. Then the corrected numerical solution is 
\begin{align}
\bq^{[1]}(s) = \bq^{[0]}(s) + \bm{\be}^{[0]}(s).
\label{}
\end{align}
Repeating the above process, one can have $\bq^{[2]},\dots,\bq^{[J]}$,
$J$ is the pre-determined number of deferred corrections. 
The convergent order of deferred correction solution along the
contour parameter is 
\begin{align}
\|\bq(s) - \bq^{[J]}(s)\| = O( (\delta s)^{m(J+1)})
\label{}
\end{align}
where $\delta s = \max\{\delta s_n\}_{n=0}^{N_s-1}$,
$m$ is the order of the chosen numerical scheme to solve
Eqns.\,\eqref{eq:matrix:form} and \eqref{eq:erreq}. For the
CN scheme, $m=2$.

In summary, under an appropriate regularity hypothesis\,\cite{lfal2013vem,Dutt2000sdc}, 
one can prove the estimator for the numerical solution $q_{\delta s, h}$, 
\begin{align}
\|q_e - q_{\delta s, h}\| = O((\delta s)^{m(J+1)} + h^{k+1})
\label{}
\end{align}
$q_e$ is the true solution of propagator, and  
$h=\max\limits_{E\in\Omega_h}\mbox{diam}\{E\}$.

\section{Numerical results}
\label{sec:rslt}

In the following numerical examples, we use linear ($k=1$) and
quadratic ($k=2$) VEMs to discretize the spatial variable. Due to
the limitation of spatial discretization order, in the time direction, we correct
the initial numerical solution one time in the SDC scheme.  
All the numerical examples are implemented based on the FEALPy package\,\cite{fealpy2021}. Halfedge data structure has also been  
integrated into FEALPy package.

\subsection{VEM with uniform mesh}



\subsubsection{Solving a parabolic equation}
\label{subsubsec:effHeat}

Solving the PDE of parabolic type is the most time-consuming part of SCFT
simulations. In this subsection, we examine the precision of our proposed method
in solving a parabolic equation.  We consider the following parabolic equation
\eqref{eq:heat}
\begin{equation}\label{eq:heat}
\left\{\begin{array}{ll}{\dfrac{\partial }{\partial s}u(x,y,s)= \dfrac{1}{2}\Delta
	u(x,y,s)},&{(x,y)\in\Omega=[0,2\pi]^2,~ s\in[0,S],}
\vspace{0.2cm}
\\ 
\dfrac{\partial}{\partial \bn} u(x,y,s)=0, & (x,y) \in \partial\Omega,
\vspace{0.2cm}
\\
{u(x,y,0)=\cos x \cos y,} 
\end{array}\right.
\end{equation}
with exact solution $u^{e}(x,y,s) = e^{-s} \cos x \cos y$.

First, we verify the convergent order of the linear and quadratic VEMs.
The CN scheme with $\delta s=1\times 10^{-4}$ is used to guarantee
enough time discretization accuracy. Tab.\,\ref{tab:spaceOrder} gives the
error and convergent order of VEM which is consistent with theoretical
results. 
\begin{table}[H]
	\centering
	\caption{The error order of VEM.} 
	\begin{tabular}{ccccc}
		\hline
		\multirow{2}*{Nodes}&\multicolumn{2}{c}{Linear VEM}&
		\multicolumn{2}{c}{Quadratic VEM}\\
		\cmidrule(lr){2-5}
		&$\|u^e(\cdot, 1) - u_h(\cdot, 1)\|_{\Omega}$& order &$\|u^e(\cdot, 1) - u_h(\cdot, 1)\|_{\Omega}$&order\\
		\hline
		289&4.7737e-02&--&1.2062e-03& --\\
		1089&1.3267e-02&1.84&1.5084e-04&2.99\\
		4225&3.4013e-03&1.96&1.8863e-05&2.99\\
		16641 &8.5563e-04&1.99&2.3582e-06&3.00\\
		\hline
	\end{tabular}
	\label{tab:spaceOrder}
\end{table}

Second, we verify the error order of the CN and SDC schemes for solving
\eqref{eq:heat}.  For the SDC scheme, we obtain a new solution $u^{[1]}$ by
correcting the initial numerical solution $u^{[0]}$ calculated by the CN scheme
just once.  For the spatial direction, we use the quadratic VEM with
$66049$ nodes to guarantee the spatial discretization accuracy.
Tab.\,\ref{tab:timeOrder} gives the convergent order of the time discretization
schemes which are also consistent with theoretical results. Note that the
error showed above is the $L^2$ error between the true solution and
numerical solution at $S=1$.
\begin{table}[H]
	\centering
	\caption{The error order of the time discretization schemes.  }
	\begin{tabular}{ccccc}
		\hline
		\multirow{2}*{$N_s$}&\multicolumn{2}{c}{CN}&
		\multicolumn{2}{c}{SDC} \\
		\cmidrule(lr){2-5}
		&$\|u^e(\cdot, s) - u_h(\cdot, s)\|_{\Omega}$&order&$\|u^e(\cdot, s) -
		u_h(\cdot, s)\|_{\Omega}$&order\\
		\hline
		4&6.0605e-03&-- &5.7514e-04& --\\                                       
		8&1.5074e-03&2.00&1.0163e-05&5.82\\                                     
		16&3.7637e-04&2.00&6.4626e-07&3.97\\                                    
		32&9.4065e-05&2.00&4.2283e-08&3.94\\
		\hline
	\end{tabular}
	\label{tab:timeOrder}
\end{table}

Third, we verify the integral accuracy of the numerical solution along with the contour variable $s$ which is required in 
solving PDEs and evaluating density functions. We use the quadratic VEM (66049
nodes) to discretize the parabolic equation \eqref{eq:heat} and obtain a
semi-discrete matrix system. Correspondingly, the exact solution of \eqref{eq:heat} can
be discretized into $u^e_h$. Then we solve the semi-discrete system using
the CN and the SDC schemes for $s\in[0,1]$ to obtain the numerical solutions
$u^{CN}_{h}$ and $u^{SDC}_h$, respectively. We integrate $u^{CN}_{h}$ and
$u^{SDC}_h$ along $s$ from $0$ to $1$ using a modified fourth-order integral
scheme \cite{press1992numerical} and the spectral integral method as
discussed in the Appendix, respectively. The integrated values are denoted
by $U_h^{CN}$ and $U_h^{SDC}$. The exact integral of $u_h^e$ along $s$ from
$0$ to $1$ can be obtained as $U^e_h$. The error is defined as 
\begin{align}
e^{M} = \| U^e_h - U^{M}_h\|_{\Omega_h},
\label{}
\end{align}
where $M\in\{CN, SDC\}$. As Tab.\,\ref{tab:intOrder} presents,
one can find that $e^{SDC}$ achieves the error level about $4\times
10^{-6}$ only requiring $8$ contour discretized nodes, while $e^{CN}$ requiring
$256$ nodes. The error value of SDC method can only be reduced to about $4
\times 10^{-6}$ due to the limitation of spatial discretization precision.
\begin{table}[H]
	\centering
	\caption{\label{tab:intOrder}
		The time integral error between the 4-order integral scheme and
		the spectral integral method}
	\begin{tabular}{ccc}
		\hline
		$N_s$&$e^{CN}$&$e^{SDC}$\\
		\hline
		4&4.2343e-03&2.4219e-04\\
		8&1.0728e-03&4.4032e-06\\
		16&2.6970e-04&4.0527e-06\\
		32&6.7675e-05&4.0519e-06\\
		64&1.7401e-05&4.0520e-06\\
		128&5.8778e-06&4.0520e-06\\
		256&4.1963e-06&4.0520e-06\\
		\hline
	\end{tabular}
\end{table}

\subsubsection{Efficiency of SCFT calculations}
\label{subsubsec:effscft}

\begin{figure}[H]
	\centering	
	\includegraphics[width=0.7\linewidth]{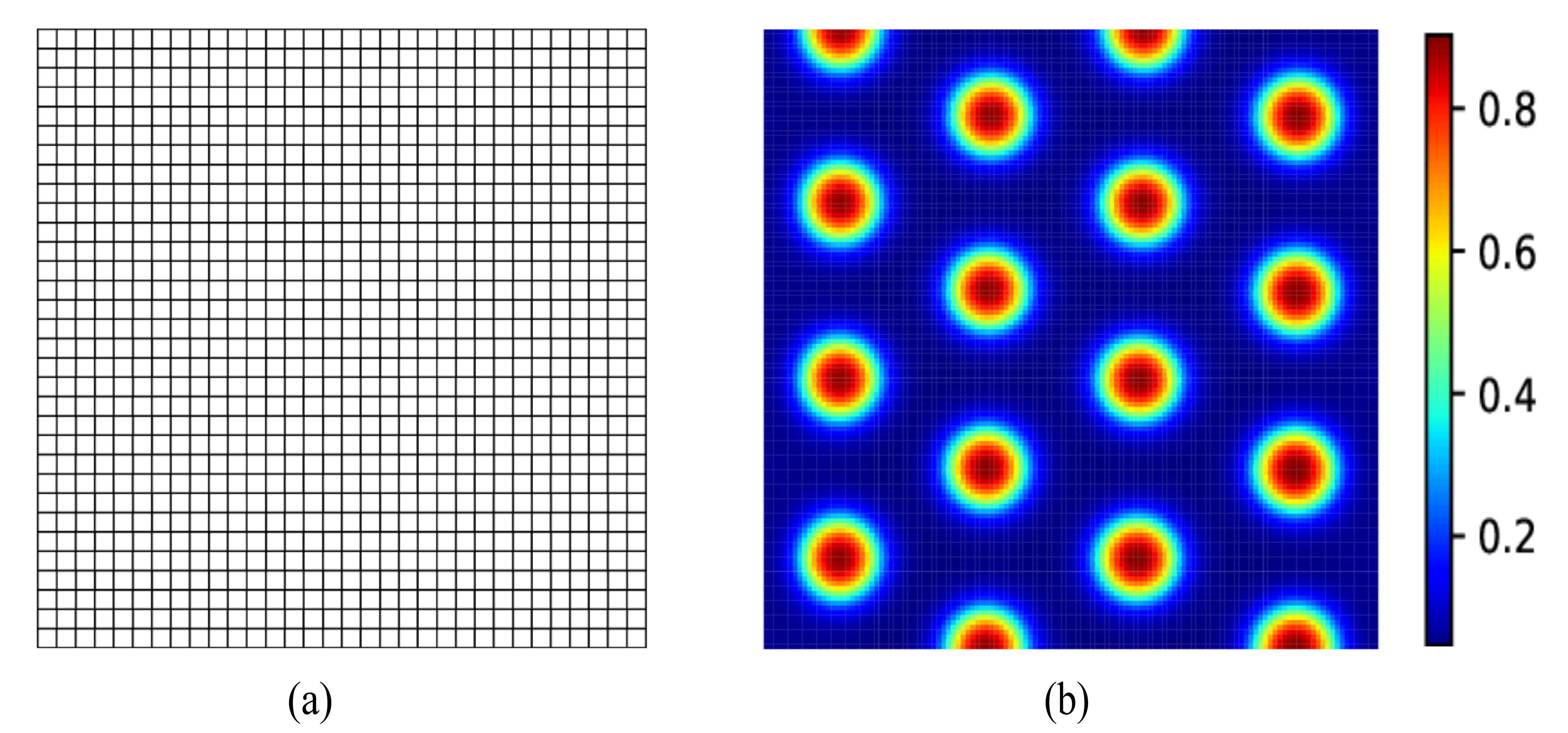}
	\caption{Cylindrical phase (b) calculated by VEM with uniform
		grid (a) when $\chi N =25$, $f=0.2$. Red colors correspond to large A-segment fractions.
	}
	\label{fig:scfteff:hex}
\end{figure}

\begin{figure}[H]
	\centering
	\setlength{\abovecaptionskip}{0.cm}
	\setlength{\belowcaptionskip}{-0.cm}
	\begin{minipage}[!htbp]{0.4\linewidth}
		\includegraphics[width=6.5cm]{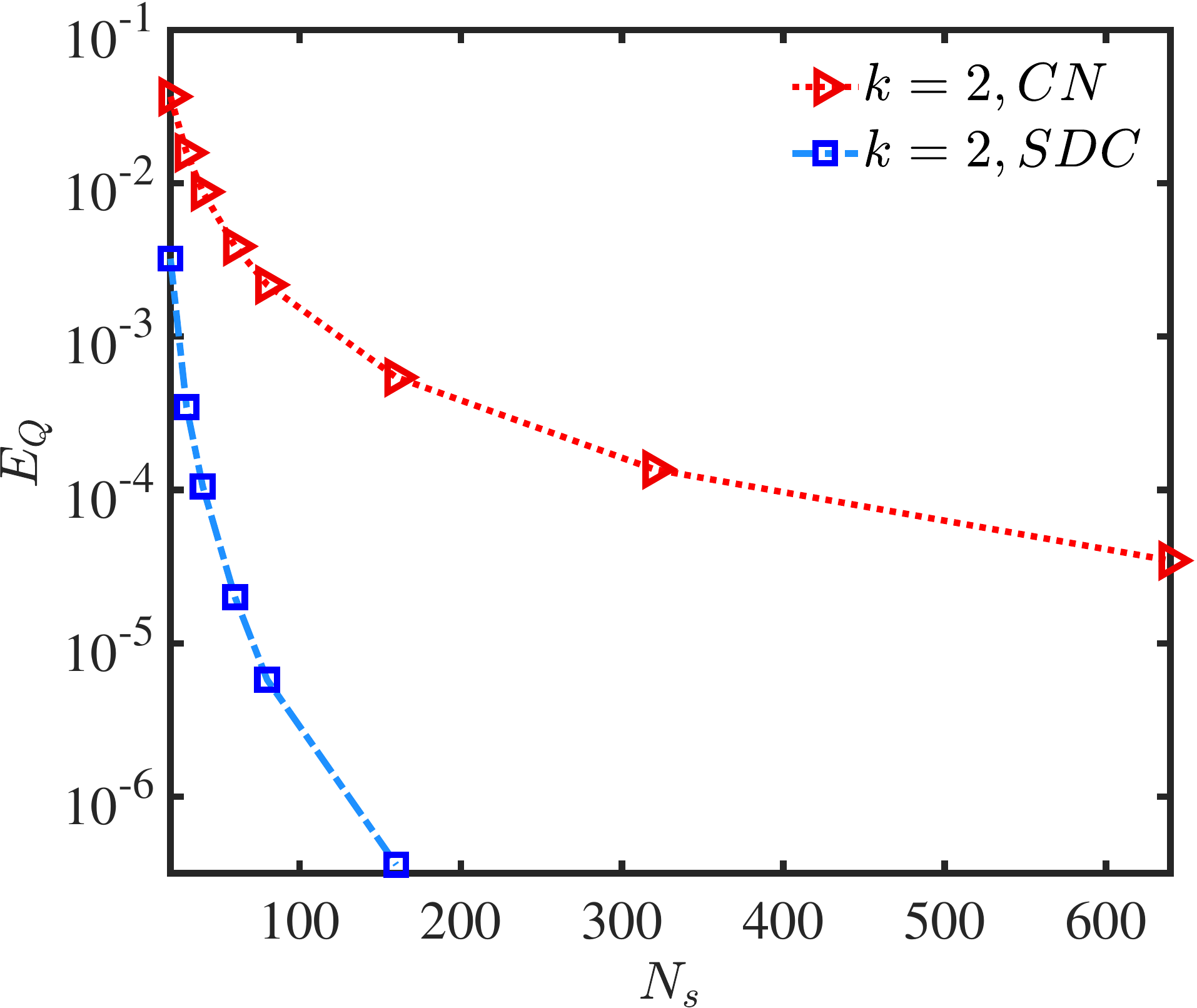}
		\caption*{(a)}
	\end{minipage}
	\hspace{0.6in}
	\begin{minipage}[!htbp]{0.4\linewidth}
		\includegraphics[width=6.5cm]{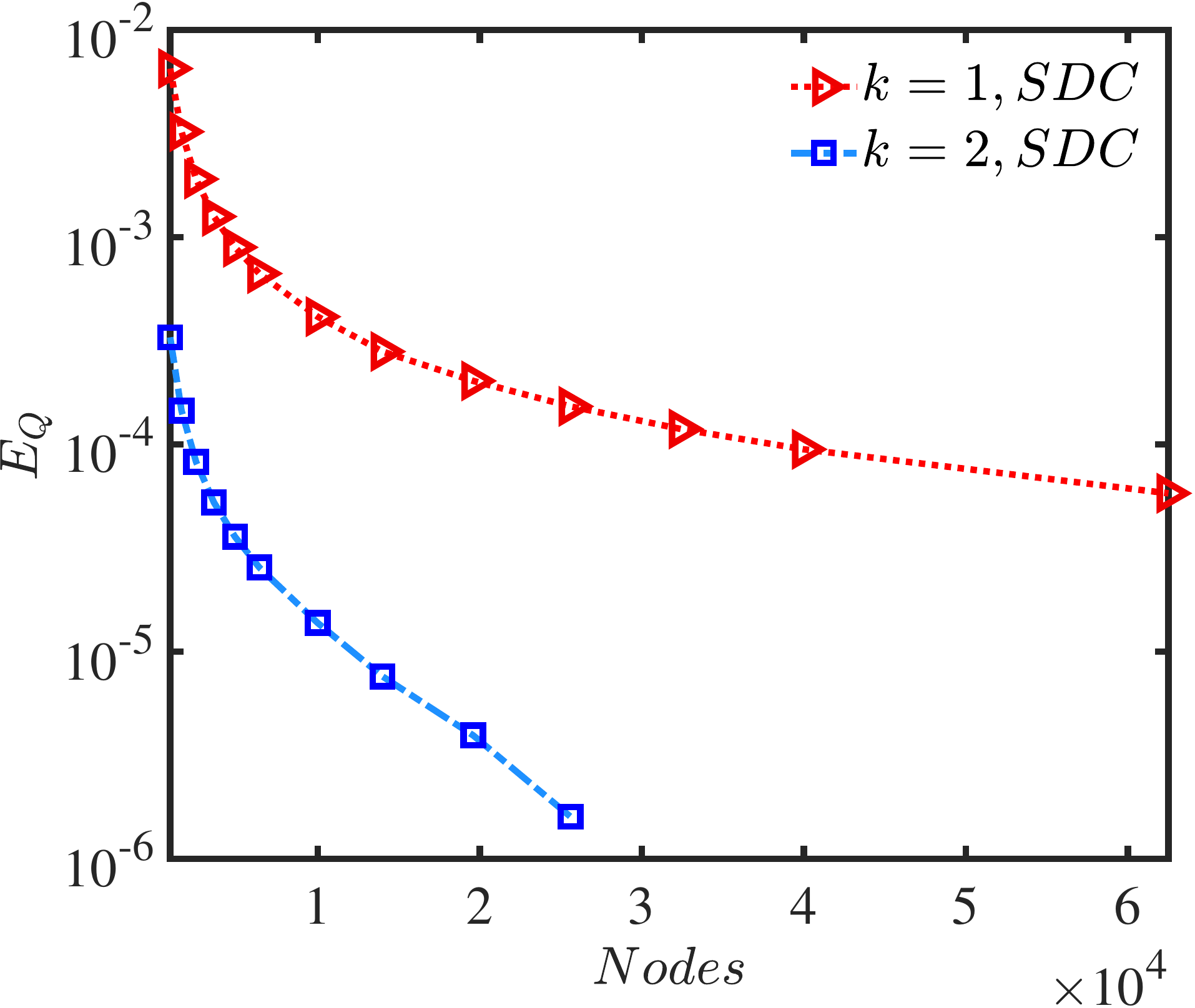}
		\caption*{(b)}
	\end{minipage}
	\caption{\label{fig:scfteff}The convergent behavior of single chain partition function $Q$ obtained by different schemes. 
		$E_Q =(Q-Q_{ref})/Q_{ref}$ is the relative error. 
		$Q_{ref}$ is the numerically exact solution. (See text for the
		details about $Q_{ref}$).
		(a) shows $Q$ obtained by the CN and SDC schemes as the contour
		points $N_s$ increase. 
		Quadratic VEM with $32400$ nodes is employed to discretize the spatial variable.
		(b) presents $Q$ computed by the linear and quadratic VEMs with an
		increase of spatial discretization points. 
		SDC scheme with $160$ points is used to discretize the contour variable.
	}
\end{figure}

To further demonstrate the performance of our proposed approach, we apply the
numerical schemes to SCFT calculations. To compare results, we need a metric
for accuracy that can be readily compared across different calculations. We use
the value of single chain partition function $Q$ as the solver's accuracy metric. Since it integrates the result of 
propagator solution, it is a measure of the entire solution. As a basis
for comparison, we use a square with an edge length of $12 R_g$ as the
computational domain. The volume fraction of $A$ is $f=0.2$, and the
interaction parameter $\chi N=25$. The computation is carried out using a
quadrilateral mesh (see Fig.\,\ref{fig:scfteff:hex}\,(a)). Correspondingly, the
convergent morphology is a cylindrical structure, as shown in
Fig.\,\ref{fig:scfteff:hex}\,(b).
First, we look at the contour discretization schemes. The goal is to have the
fewest number of contour points necessary for a desired accuracy. The
quadratic VEM with $32400$ nodes is used to guarantee enough spatial
discretization accuracy. $Q_{ref}$ in Fig.\,\ref{fig:scfteff}(a) is numerically
obtained by the SDC scheme with $160$ contour points. Fig.\,\ref{fig:scfteff}(a)
shows the convergent information  of $Q$ for the CN and SDC schemes, as discussed above. The
SDC method converges faster than the CN scheme to a prescribed precision. 

Second, we observe the numerical behavior of linear and quadratic VEMs in the
SCFT simulation. From the above numerical tests (see
Fig.\,\ref{fig:scfteff}(a)), one can see that using the SDC scheme with $160$
discretization points can guarantee enough accuracy in the contour direction.  So in
the following computations, we use a high-precision numerical $Q_{ref}$ as the
exact value, which is obtained by the quadratic VEM with $32400$ nodes and SDC
scheme with $160$ points.  Fig.\,\ref{fig:scfteff} (b) shows the $Q$ values with
different spatial discretization points of linear and quadratic VEMs. It is easy
to see that the quadratic VEM is more accurate than the linear VEM as theory
predicts.  Therefore, in the following calculations, we always adopt the
quadratic VEM and the SDC scheme.

\subsubsection{General domains with general polygonal mesh}
\begin{figure}[H]
	\centering
	\includegraphics[width=1.02\linewidth]{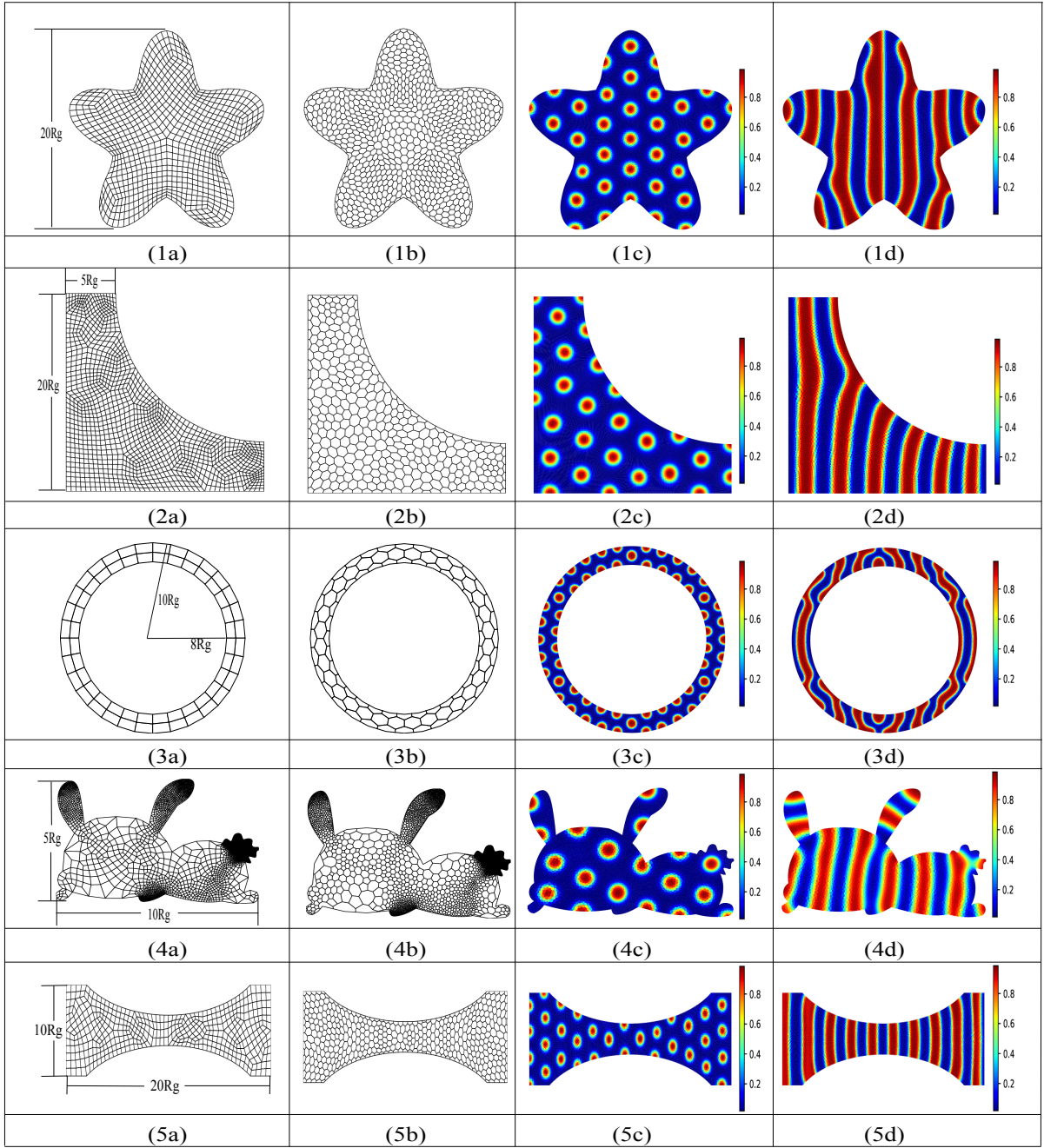}
	\caption{\label{fig:general} 
		The self-assembled patterns in general domains through SCFT simulation, including
		(1). Flower shaped plane; (2). Curved-L shaped plane; (3). Ring domain; (4).
		Rabbit-shaped plane; and (5). Dumbbell plane.
		Red colors correspond to A-segment fractions. The first and second columns
		present the schematic mesh of quadrangular and polygonal meshes, respectively. 
		The simulating diblock copolymer systems contain (1c) $[\chi N, f]=[25, 0.2]$,
		(1d) $[\chi N, f]=[15, 0.5]$,
		(2c) $[\chi N, f]=[25, 0.2]$, (2d) $[\chi N, f]=[15, 0.5]$,
		(3c) $[\chi N, f]=[30, 0.2]$, (3d) $[\chi N, f]=[14, 0.5]$,
		(4c) $[\chi N, f]=[30, 0.3]$, (4d) $[\chi N, f]=[14, 0.5]$,
		(5c) $[\chi N, f]=[25, 0.2]$, and (5d) $[\chi N, f]=[15, 0.5]$. 
		The number of nodes of the mesh and converged Hamiltonian values can be found in Tab.\,\ref{tab:general}.
	}
\end{figure}

One advantage of VEM can use the arbitrary approximate geometry
domain with general polygonal meshes. 
Fig.\,\ref{fig:general} presents these
results on five different two-dimensional domains divided by quadrilateral and
polygonal elements, respectively.  The same convergent structure and almost the
same Hamiltonian value can be obtained for these two kinds of meshes, as
shown in Fig.\ref{fig:general} and Tab.\,\ref{tab:general}.
\begin{table}[H]
	\centering
	\caption{The number of nodes of different meshes used in SCFT calculations
		for five different domains as shown in
		Fig.\,\ref{fig:general} and corresponding converged Hamiltonian values.}
	\begin{tabular}{ccccc}
		\hline
		\multirow{2}*{Domain}&\multirow{2}*{Mesh}&\multirow{2}*{Nodes}&\multicolumn{2}{c}{Hamiltonian}\\
		\cmidrule(lr){4-5}
		&&&(c)&(d)\\
		\hline
		
		\multirow{2}*{Fig.\,\ref{fig:general}\,(1)}&(a)&13041&-2.3742&-1.7388\\
		&(b)&22560&-2.3754&-1.7398\\ 
		\hline
		\multirow{2}*{Fig.\,\ref{fig:general}\,(2)}&(a)&10720&-2.3720&-1.7440\\
		&(b)&20273&-2.3765&-1.7382\\
		\hline
		\multirow{2}*{Fig.\,\ref{fig:general}\,(3)}&(a)&7014&-3.1410&-0.1874\\
		&(b)&6510&-3.1409&-0.1873\\
		\hline
		\multirow{2}*{Fig.\,\ref{fig:general}\,(4)}&(a)&30182&-3.1440&-0.1900\\
		&(b)&34587&-3.1448&-0.1901\\
		\hline
		\multirow{2}*{Fig.\,\ref{fig:general}\,(5)}&(a)&7601&-2.3670&-1.6797\\
		&(b)&13824&-2.3718&-1.6883\\
		\hline
	\end{tabular}
	\label{tab:general}
\end{table}

\subsection{VEM with adaptive mesh}
\label{subsec:rltsAdaptiv}

In this subsection, we will demonstrate the efficiency of adaptive VEM from three
parts: 1) the less computational cost to obtain prescribed accuracy compared with uniform mesh; 2)
the application to strong segregation systems; 3) general domains with adaptive polygonal mesh. As discussed in
Sec.~\ref{subsubsec:effHeat}, the quadratic VEM is more accurate than the linear
one. Therefore, only quadratic VEM is used in the adaptive process.
Meanwhile, the SDC scheme is chosen to discretize
the contour variable with $100$ points.  

\subsubsection{Efficiency}

\begin{figure}[H]
	\centering
	\includegraphics[width=0.7\linewidth]{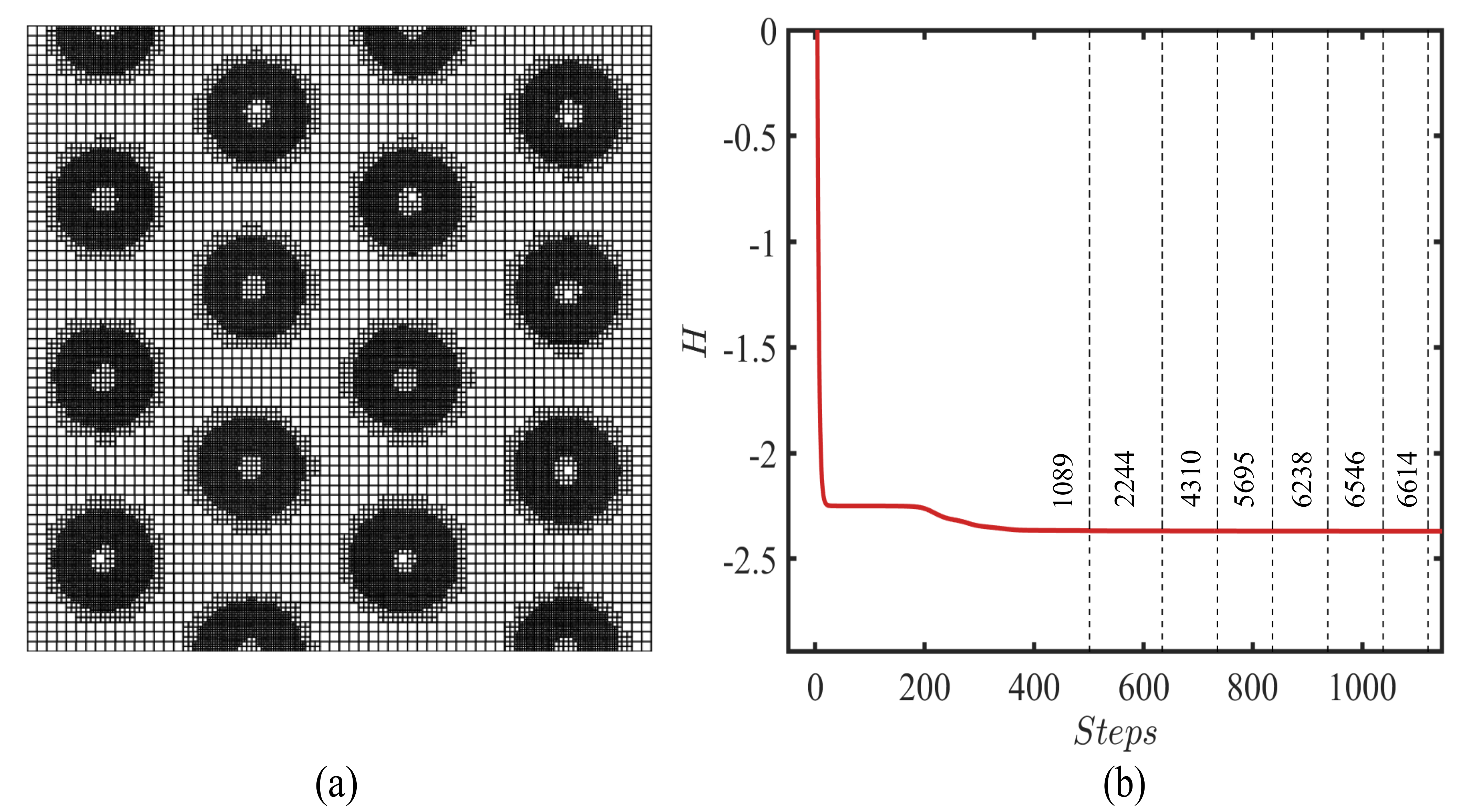}
	\caption{\label{fig:avem:eff} 
		(a) The converged adaptive mesh. (b) The numerical behavior of Hamiltonian $H$.
		The numbers between two dotted lines represent the number of
		spatial nodes in the adaptive process.}
\end{figure}

First, we take $\chi N = 25$ and $f=0.2$ as an example to demonstrate the
efficiency of adaptive VEM. The computational domain is a square with
an edge length of $12 R_g$. The square domain's uniform mesh with $1089$ nodes is used
to model the system at the start stage. Then adaptive VEM is launched
when the iteration reaches the maximum step $500$ or the reference value of the estimator
$\eta_{ref}<0.1$, where
$$
\eta_{ref} = \sigma(\eta_E)/(\max(\eta_E)-\min(\eta_E)).
$$
$\sigma(\eta_E)$ is the standard deviation of $\eta_E$, estimator $\eta_E$ see
Eqn.\,\eqref{eq:avem:error}). The adaptive process will be terminated when the
the successive Hamiltonian difference is smaller than $1.0\times 10^{-6}$.
Fig.\,\ref{fig:avem:eff} (a) gives the final adaptive mesh which includes
$6684$ nodes.  Fig.\,\ref{fig:avem:eff} (b) shows the convergent tendency of
Hamiltonian $H$ of the adaptive process.  The finally converged morphology has been
shown in Fig.\,\ref{fig:scfteff:hex}\,(b).  It can be seen that the Hamiltonian
value efficiently converges by the cascadic adaptive method and refined meshes
concentrate on the shape interface. 

\begin{figure}[H]
	\centering
	\setlength{\abovecaptionskip}{0.cm}
	\setlength{\belowcaptionskip}{-0.cm}
	\begin{minipage}[!htbp]{0.4\linewidth}
		\includegraphics[width=6.5cm]{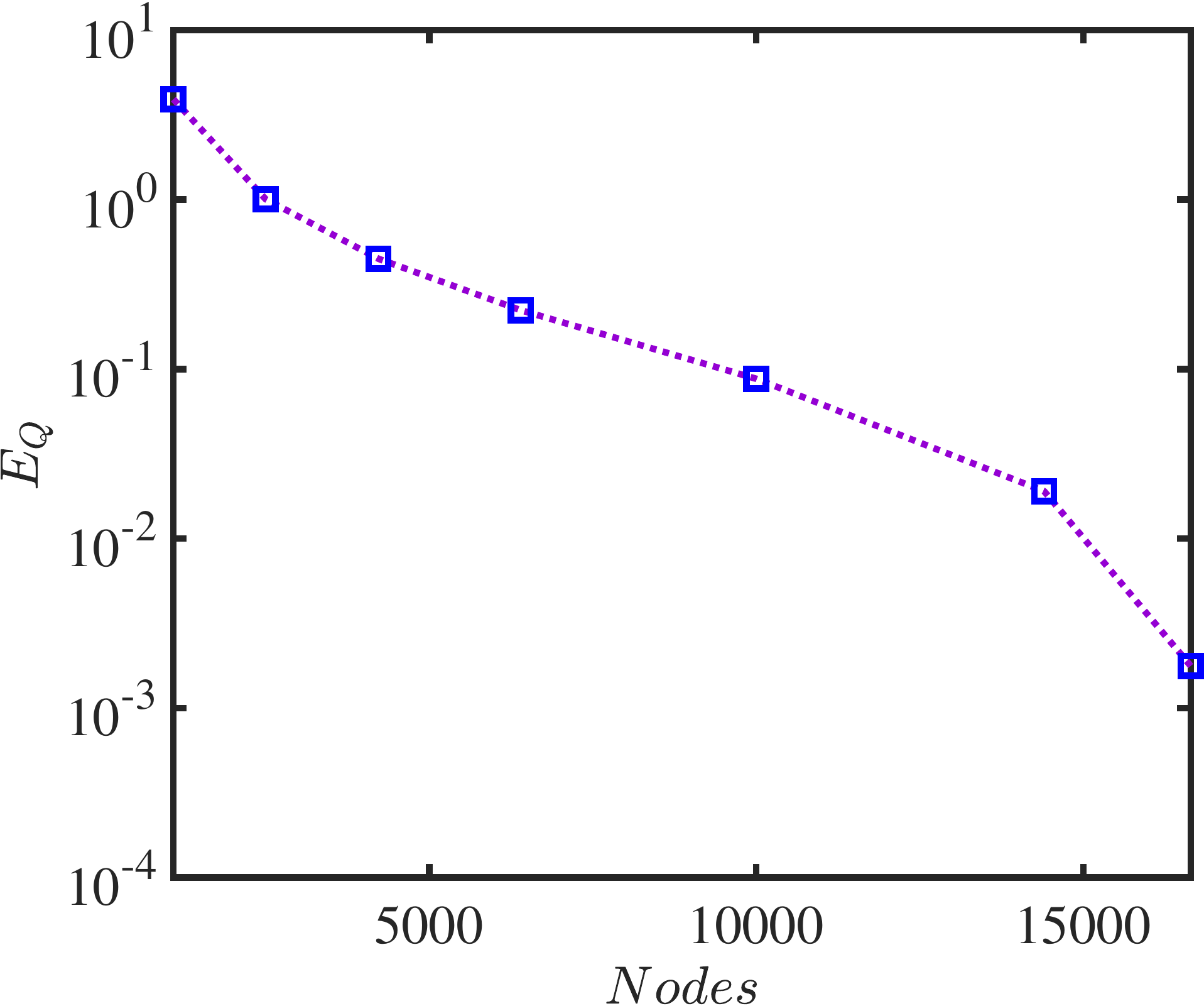}
		\caption*{(a)}
	\end{minipage}
	\hspace{0.6in}
	\begin{minipage}[!htbp]{0.4\linewidth}
		\includegraphics[width=6.5cm]{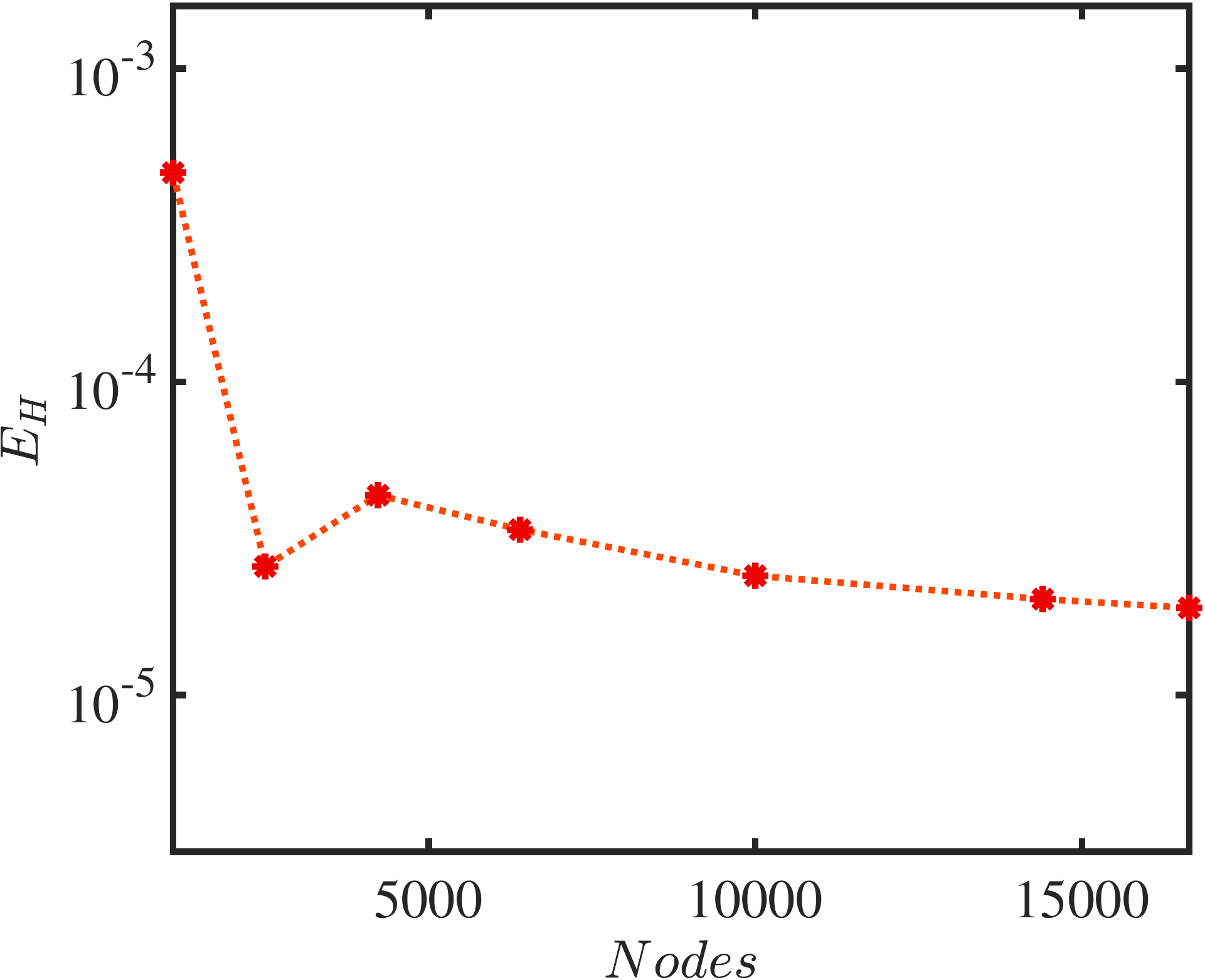}
		\caption*{(b)}
	\end{minipage}
	\caption{\label{fig:compare} 
		The convergence results of VEM with adaptive and uniform mesh
		when $\chi N =25, f=0.2$.  The differences of (a) single partition function
		$Q$, $E_Q =(Q-Q_{adap})/Q_{adap}$ and of (b) Hamiltonian value$H$, $E_H =
		(H-H_{adap})/H_{adap}$. $Q$ and $H$ are obtained with uniform mesh, while
		$Q_{adap}$ and $H_{adap}$ are calculated with the adaptive approach.
	}
\end{figure}

We also compared the simulation results of VEM with adaptive and uniform mesh. 
Fig.\,\ref{fig:compare} shows the numerical behaviors of single chain
partition function $Q$ and Hamiltonian $H$ as the nodes increase.
Tab.\,\ref{tab:avem-uniform} gives the corresponding converged values of $Q$ and
$H$.
From these results, one can find that the uniform mesh's results indeed gradually
converge to that of adaptive VEM. However, there exists a small gap between the
results of the two methods. The reason is that the adaptive VEM puts more meshes on the sharp interface and obtained a relatively accurate solution.
The minimum element size of the
adaptive mesh in the above calculation is $h_{min}=0.0469Rg$. 
While the uniform mesh method with the same element size $h_{min}$ 
requires about $65000$ nodes, which is about ten times the adaptive approach.

\begin{table}[H]                                                                 
	\centering
	\caption{\label{tab:avem-uniform} 
		The convergence values of $Q$ and $H$ obtained by VEM with adaptive and uniform meshes.}
	\begin{tabular}{cccc}                                                
		\hline                                                                  
		Mesh &Nodes & Q & H\\
		\hline
		Adaptive &  6684   & 4.2295e+02 & -2.369403\\
		Uniform  &   16641 & 4.2373e+02 & -2.369448\\
		\hline
	\end{tabular}
\end{table}

\subsubsection{Strong segregation systems}

Next, we apply the adaptive VEM to simulate strong segregation systems,
i.e., large interaction parameter $\chi N$, also in the square domain with an edge
length of $12R_g$. For the strong segregation case, the interface thickness becomes
narrower.  Therefore the adaptive method is more suitable than the uniform mesh approach to
catch these narrower interfaces. When simulating the strong segregation system, the
initial values are obtained by the relatively weak
segregation system's converged results. Tab.\,\ref{tab:avem-adaptive} presents the numerical results
of $\chi N$ from $25$ to $60$ and $f=0.2$. From these results, one can find the
advantages of the adaptive VEM as $\chi N$ increases, including a mild
increase of mesh nodes and fewer iteration steps.

\begin{table}[H]                                                                
	\centering
	\caption{\label{tab:avem-adaptive} 
		Numerical results by the adaptive VEM for strong segregation systems.}
	\begin{tabular}{ccccc}                                                
		\hline
		$\chi N$ & Step & Nodes & H\\\hline
		25 & 1146 & 6684 &  -2.369403\\
		30 & 78  & 9037  &  -3.149607\\
		35 & 89  & 13443 &  -4.020791\\
		40 & 74  & 17649 &  -4.946249\\
		45 & 75  & 19741 &  -5.907039\\
		50 & 75  & 20480 &  -6.892386\\
		55 & 73  & 20641 &  -7.895548\\
		60 & 61  & 20690 &  -8.911902\\
		\hline
	\end{tabular}
\end{table}

\begin{figure}[H]
	\centering
	\includegraphics[width=1.03\linewidth]{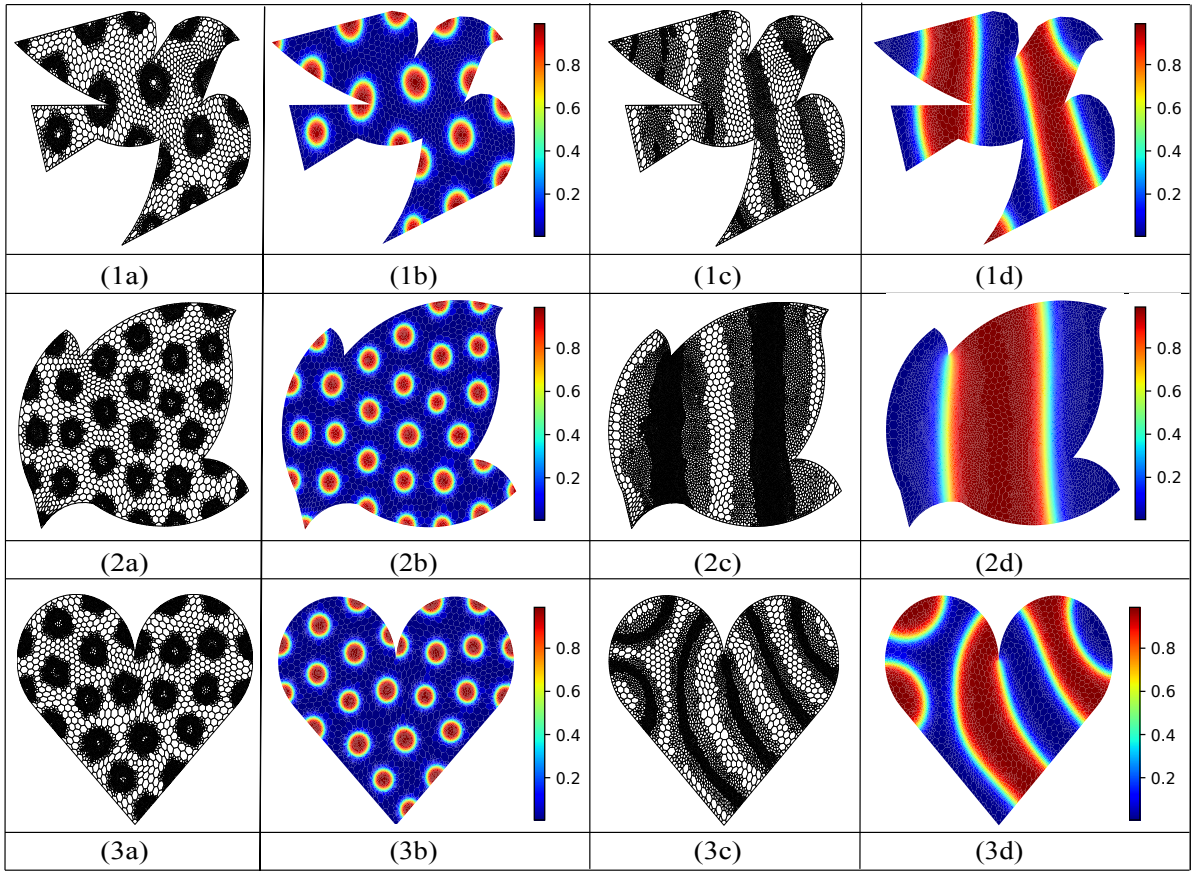}
	\caption{\label{fig:other} 
		The self-assembled patterns of strong segregation systems obtained by the
		adaptive VEM, (a)(c) adaptive meshes, (b) spotted phases, (d) lamellar
		phases. The model parameters are
		(1b) $[\chi N, f]=[35, 0.2]$, (1d) $[\chi N, f]=[30, 0.5]$,
		(2b) $[\chi N, f]=[40, 0.2]$, and (2d) $[\chi N, f]=[40, 0.5]$.
		(3b) $[\chi N, f]=[40, 0.2]$, and (2d) $[\chi N, f]=[30, 0.5]$.
		Red colors correspond to large A-segment fractions.
	}
\end{figure}

Finally, we apply the adaptive VEM to the strong segregation
systems on more complicated domains,  including two kinds of structures, spotted
phases when $f =0.2$ and lamellar phases when $f =0.5$.
Fig.\,\ref{fig:other} presents the adaptive meshes and converged morphologies.
The corresponding nodes of uniform mesh are estimated by the minimum mesh size of the adaptive mesh.  
A comparison demonstrates that the adaptive method can greatly reduce  the number of nodes as shown in Tab.\,\ref{tab:other}.



\begin{table}[H]                                                                
	\centering
	\caption{\label{tab:other} The corresponding data on three different planes
		as shown in Fig.\,\ref{fig:other}.
		$Node_{adap}$ and $h_{min}$ are the number of nodes and minimum grid size of the 
		adaptive mesh, respectively. $Node_{uni}$ is the number of nodes of
		the uniform mesh with the same cell size estimated by $h_{min}$. $C_{save} =
		1- Node_{adap}/Node_{uni} $. }
	\begin{tabular}{ccccc}                                                
		\hline
		Mesh & $h_{min}$ & $Node_{adap}$ & $Node_{uni}$ &
		$C_{save}(\%)$ \\
		\hline
		Fig.\,\ref{fig:other}(1)(a)  & 3.54e-02 & 8154 & 12956 &  37.1\%\\
		Fig.\,\ref{fig:other}(1)(c)  & 1.17e-02 & 7591  & 26727   &  71.6\%\\
		Fig.\,\ref{fig:other}(2)(a)  & 2.98e-02 & 19118 & 30169  &  36.6\%\\
		Fig.\,\ref{fig:other}(2)(c)  & 2.12e-03 & 24138 & 116533  &  79.3\%\\
		Fig.\,\ref{fig:other}(3)(a)  & 3.04e-02 & 14914 & 27729  &  46.2\%\\
		Fig.\,\ref{fig:other}(3)(c)  & 1.33e-02 & 8399  & 25763  &  67.4\%\\
		\hline
	\end{tabular}
\end{table}


\section{Conclusion}
\label{sec:conclusion}
In this paper, we propose an efficient numerical method to solve the polymer
SCFT model on arbitrary domains based on the VEM.  We have developed an adaptive
method equipped with a new $Log$ marking strategy that can make full use of the
information of numerical results and save the SCFT iterations significantly. Using the
halfedge data structure, we can apply an adaptive method to refine and coarsen
arbitrary polygonal grids.  
The SDC method is also used to discretize the contour variable. 
The resulting method can obtain a high-accuracy numerical solution with fewer spatial and contour nodes.
Numerical results demonstrate that the adaptive VEM even saves the computational amount up to $79.3\%$ in solving a strong segregation lamellar system
compared with the uniform mesh method. 
In this work, we have applied our algorithms to two-dimensional SCFT calculations. In future work, we aim to develop the adaptive VEM method to investigate three-dimensional SCFT problems.

\section*{Appendix: Spectral integral method along the contour variable $s$}
\begin{appendix}
\label{subsec:spint}

In this Appendix, we discuss the Chebyshev-node interpolatory quadrature method
\cite{ceniceros2019efficient} 
to integrate the residual error $\bm\gamma^{[0]}(s)$ of
Eqn.\,\eqref{eq:residual} for the contour variable $s$, which has
spectral accuracy for smooth integrand\,\cite{trefethen2008guass}. 
The proposed scheme can also be applied to evaluating the density operators \eqref{eq:phiA} and
\eqref{eq:phiB}. These problems can be summarized to the
following integral
\begin{equation}
    \int_{a}^{b} g(s)\, ds,
\end{equation}
where the integrand $g(s)$ is a smooth function.
After changing variables, the general integral becomes
\begin{equation}\label{Eqn:theta}
    \int_{a}^{b} g(s) \,ds 
     = \frac{b-a}{2}\int_{0}^{\pi} g(-\cos\theta)\sin\theta \,d\theta,
\end{equation}
where $\theta \in [0,\pi]$. 
The interpolate polynomial of $g$ at 
Chebyshev node $\theta_j=j\pi/N_s$, $j=0,1,\dots,N_s$.
\begin{equation}
    g(-\cos \theta_j) \approx \frac{a_{0}}{2}+\sum_{k=1}^{N_s-1}
	a_{k} \cos (k \theta_j) + \frac{1}{2} a_{N_s} \cos (N_s
	\theta_j),
\end{equation}
where
\begin{equation}
     a_{k}=\frac{1}{N_s}g\left(-\cos\theta_{0}\right)\cos (k
 \theta_{0})+\frac{2}{N_s}
	{\sum_{j=1}^{N_s-1}}g\left(-\cos\theta_{j}\right)\cos (k
    \theta_{j})+\frac{1}{N_s}g\left(-\cos\theta_{N_{s}}\right)\cos\left(k\theta_{N_s}
\right),
\end{equation}
$k=0, \ldots, N_s$. In practice, coefficients $a_0,a_1,\ldots,a_{N_s}$ are calculated
by the fast discrete cosine transform.
\begin{equation}
\int_{0}^{\pi} g(-\cos\theta) \sin \theta \, d \theta=\frac{a_{0}}{2}
\int_{0}^{\pi} \sin \theta \,d \theta+\sum_{k=1}^{N_s-1} a_{k} \int_{0}^{\pi} \cos
k \theta \sin \theta \, d \theta+\frac{a_{N_s}}{2} \int_{0}^{\pi}
\cos N_s \theta \sin \theta \,d \theta.
\end{equation}
Due to $\cos k \theta \sin \theta=\dfrac{1}{2}[\sin (1+k) \theta+\sin (1-k)
\theta]$, we have
\begin{equation}
	\int_{a}^{b} g(s)\,ds = \frac{b-a}{2}\int_{-1}^{1}g(t)\,dt
	\approx
    \left\{\begin{array}{ll}\dfrac{b-a}{2}\left[ a_0 +
		\sum\limits_{k=2 \atop k \text { even }}^{{N_s}-2}
            \dfrac{2a_{k}}{1-k^{2}}+\dfrac{a_{N_s}}{1-N_s^{2}}\right],
			& \quad N_s  \ \text{is even},\\
        \\
    \dfrac{b-a}{2}\left[ a_0 + \sum\limits_{k=2 \atop k \text { even
	}}^{{N_s}-1}\dfrac{2 a_{k}}{1-k^{2}}\right], & \quad N_s \ \text{is odd}.\\
\end{array}\right.
\end{equation}

\end{appendix}

\end{document}